\documentclass[journal,twoside]{IEEEtran}
\newif\ifmtpro
\mtprofalse
\ifmtpro
 \usepackage[subscriptcorrection,slantedGreek,mtpbb,mtpfrak,zswash]{mtpro2}
 
\else
 \usepackage[italic]{mathastext}
 \usepackage{amssymb}
 \newcommand{\coloneq}{:=}
 \newcommand{\eqcolon}{=:}
 \newcommand{\ccases}[1]{\begin{cases} #1 \end{cases}}
\fi

\usepackage{amsmath,cite}

\usepackage{graphicx,color,psfrag}
\usepackage[tight,scriptsize]{subfigure}
\usepackage{arydshln}
\newtheorem{theorem}{Theorem}[section]
\newtheorem{lemma}[theorem]{Lemma}
\newtheorem{proposition}[theorem]{Proposition}
\newtheorem{corollary}[theorem]{Corollary}
\newtheorem{remark}{Remark}[section]

\DeclareMathOperator\trace{tr\hspace*{.5pt}}

\DeclareMathOperator\sign{sign}
\DeclareMathOperator\real{Re}
\DeclareMathOperator\im{Im}

\newcommand\norm[1]{\ensuremath{\lVert#1\rVert}}
\newcommand\normf[1]{\ensuremath{\lVert#1\rVert_\textsc f}}
\newcommand{\inpr}[2]{\ensuremath{\langle#1,#2\rangle}}

\newcommand\D{\mathrm{d}}
\newcommand\E{\mathrm{e}}
\newcommand\J{\hspace{.5pt}\mathrm{j}}
\newcommand\Htwo{\ensuremath{H^2}}
\newcommand\Hinf{\ensuremath{H^\infty}}

\newcommand\mmatrix[2][cccccc]{\ensuremath{\left[\begin{array}{#1}#2\end{array}\right]}}

\newcommand\LFTl[2]{\ensuremath{\mathcal{F}_{\text l}(#1,#2)}}
\newcommand\LFTlL[2]{\ensuremath{\mathcal{F}_{\text l}\left(#1,#2\right)}}
\newcommand\LFTu[2]{\ensuremath{\mathcal{F}_{\text u}(#1,#2)}}

\newcommand\vd[1]{\ensuremath{\bar{#1}}}     
\newcommand\vl[1]{\ensuremath{\breve{#1}}}   

\makeatletter
\xdef\@endgadget#1{\unskip\nobreak\hfil\penalty50\hskip1em\hbox{}\nobreak\hfil#1\parfillskip=0pt\finalhyphendemerits=0\par}
\def\@finsymbol{$\circ$}
\def\fin{\@endgadget{\@finsymbol}}
\def\@Endofsymbol{$\triangledown$}
\def\Endoftheorem{\@endgadget{\@Endofsymbol}}
\makeatother

\title{Intermittent Redesign of Analog Controllers via the Youla Parameter}

\author{Leonid Mirkin\thanks{Faculty of Mechanical Eng., Technion---IIT, Haifa 32000, Israel. E-mail: \textsl{mirkin@technion.ac.il}.}}

\begin{document}

\maketitle

\begin{abstract}
 The paper studies digital redesign of linear time-invariant analog controllers under
 intermittent sampling. The sampling pattern is only assumed to be uniformly bounded,
 but otherwise irregular and unknown a priori. The contribution of the paper is
 twofold. First, it proposes a constructive algorithm to redesign any analog
 stabilizing controller so that the closed-loop stability is preserved.  Second, it is
 shown that when applied to (sub)\,optimal $\boldsymbol{\Htwo}$ and
 $\boldsymbol{\Hinf}$ controllers, the algorithm produces (sub)\,optimal sampled-data
 solutions under any a priori unknown sampling pattern. The proposed solutions are
 analytic, computationally simple, implementable, and transparent.  Transparency pays
 off in showing the optimality, under a fixed sampling density, of uniform sampling
 for both performance measures studied.
\end{abstract}



\section{Introduction} \label{sec:intro}

\begin{figure}[!tb]
 \centering\includegraphics[scale=0.667]{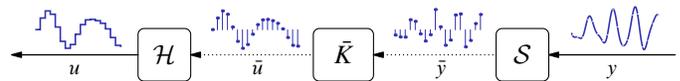}
 \caption{Generic sampled-data controller as the cascade of an A/D converter (sampler) $\mathcal S$, a pure discrete-time part $\vd K$, and a D/A converter (hold) $\mathcal H$} \label{fig:sdc}
\end{figure}

The term ``digital redesign'' refers to problems of approximating analog controllers
by sampled-data ones, i.e.\ controllers that can be realized as the cascade of a
sampler (A/D converter), a pure discrete element, and a hold (D/A converter) as shown
in Fig.\,\ref{fig:sdc}. This approach has been widely employed in designing digital
controllers for analog plants, not least because it facilitates the direct use of
analog insights in the design. The reader is referred to \cite[Ch.\:8]{AW:97} and
\cite[Ch.\:3]{ChFr:95} for expositions of ideas in the field and further references.

A common digital redesign setup is to assume a regular (say, constant) sampling rate,
fixed A/D and D/A parts (say, the ideal sampler and the zero-order hold, respectively,
as in Fig.\,\ref{fig:sdc}), and choose a discrete-time part that mimics the structure
of the analog prototype. But these choices are, to some extent, a legacy of
technological and methodological limitations of early computer-controlled systems.
Nowadays, with the advent of affordable DSP technology and a trend to distribute
information processing, the accents are changing.

First, the use of traditional A/D and D/A converters might no longer be preordained.
There may be enough local computational power to pre-process measurements and
post-process control commands. Model-based modifications of control signal during the
intersample, dubbed the generalized hold, were exploited in \cite{K:87} (in fact, an
application of a generalized hold mechanism to the digital redesign problem was
already proposed in \cite{YKS:74}), with the philosophy to circumvent limitations of
linear control. This philosophy was then criticized in \cite{FG:94}. Optimal design of
generalized sampler and hold, which are not prone to the problems presented in
\cite{FG:94}, was pioneered in \cite{Tad:92}, see also \cite{MRP:99I}. Lately, there
is a renewed interest in this subject, see e.g.\ \cite{A:08,GAM:14} and the references
therein.

Second, there have been rapidly growing activities in systems with intermittent
sampling. This is motivated by networked control systems \cite{GAM:14} and potential
advantages in employing event-based feedback \cite{A:08,HJT:12}. Although the results
might not study digital redesign explicitly (an exception is \cite{GW:10}), many of
them effectively deal with these problems. Of a special interest for us are approaches
that make use of the simulated analog closed-loop system to generate control signals
during intersample intervals of irregular lengths and, if not the whole state is
measured, adjust an analog state estimator upon arrival of new samples. This direction
is exposed in \cite{GAM:14}. See also \cite{FrGr:95} for apparently the first
appearance of such an idea in the control literature and \cite{GLLG:11} and the
references therein for its use in human control, although these two references offer
neither proofs of stability nor performance analyses. It is worth emphasizing that
many methods, which use intermittent sampling, augment the original analog controller,
so that its discretized version may be more complex. This departure from the
conventional modus operandi reflects the changing accents mentioned above: more
emphasis is placed on the information exchange between system components and the form
of A/D and D/A converters is less restrictive.

Tackling systems with intermittent sampling events might be a challenge, owing to
their time-varying nature and switches between closed- and open-loop regimes. This is
true in handling the closed-loop stability and even more so in analyzing performance.
Consequently, results are frequently either conservative or apply only to simple
dynamics. The full access to the plant state is a recurrent assumption. There appear
to be no non-conservative and transparent methods of optimal control design for
general linear problems with general sampling patterns. Besides, although the use of
unorthodox hold and sampling elements has proved useful, their structures are often
justified only empirically. The apparent qualitative difference from systems with
periodic sampling brought about different analysis tools, like continuous-time
Lyapunov methods.

One of the goals of this paper is to demonstrate that concepts and tools developed for
sampled-data problems with periodic sampling might still be powerful in addressing
stability and performance problems under intermittent sampling. It is shown that the
ideas of \cite{MR:97}, which exploit properties of conventional sampled-data systems
in the \emph{lifted domain}, extend to systems with intermittent sampling.
Specifically, \cite{MR:97} shows that the set of all causal finite-dimensional
sampled-data systems corresponds to the set of strictly causal systems in the lifted
domain. This result facilitates extracting sampled-data controllers from various
analog controller parametrizations. By extending the result to the intermittent
sampling setup, the following set of redesign problems is addressed:
\begin{enumerate}[\IEEEsetlabelwidth{}]
\item An approach to digitally redesign given analog stabilizing controllers is put
 forward. By embedding such controllers into the analog Youla parametrization setup,
 all stabilizing sampled-data controllers are characterized. This yields a systematic
 algorithm to construct a stabilizing controller under any, even unknown a priori,
 sampling pattern.
\item Intermittent redesign methods for analog \Htwo\ and \Hinf\ (sub)\,optimal
 controllers are proposed. They result in non-conservative optimal designs under no
 limitation on the sampling pattern. Performance levels attainable by the resulting
 sampled-data controllers are transparent functions of sampling times. As a result, it
 is proved that the uniform sampling is both \Htwo\ and \Hinf\ optimal among all
 sampling patterns of a given density.
\end{enumerate}
Remarkably, the offline computational complexity of the algorithms above is
independent of the sampling pattern. Also, the resulting sampler and hold are
justified performance wise. To the best of my knowledge, these are the first
non-conservative and computationally tractable results for general linear problems
with unrestricted sampling patterns.

The paper is organized as follows. After presenting some preliminary results about the
Youla parametrization and lifting in Section~\ref{sec:prelim}, the class of
sampled-data controllers is characterized in the lifted domain in
Section~\ref{sec:sdc}. This result is then used to address the stabilization problem
in Section~\ref{sec:stab}, where a parametrization of all sampled-data stabilizing
controllers for an arbitrary sampling pattern is presented (Theorem~\ref{th:Ksd}) and
some of their properties are discussed. The next section is devoted to the
performance-based discretizations, in the \Htwo\ (\S\ref{sec:H2}) and \Hinf\ 
(\S\ref{sec:Hinf}) senses. Section~\ref{sec:ex} shows how the proposed approach can be
applied to the \Hinf\ loop shaping method of \cite{McG:90} and illustrates this
procedure by a numerical example. Concluding remarks are provided in
Section~\ref{sec:concl} and the Appendix contains some more technical proofs.

\medskip
\paragraph*{Notation}

The sets of non-negative integers and reals are denoted as $\mathbb Z^+$ and $\mathbb
R^+$, respectively. The transpose of a matrix $M$ is denoted as $M'$ and, for square
matrices, $\trace(M)$ and $\rho(M)$ stand for the trace and the spectral radius of
$M$. $\LFTl\Phi\Omega$ and $\LFTu\Phi\Omega$ read as the lower and upper
linear-fractional transformations of $\Omega$ by $\Phi$, respectively, see
\cite[Ch.\,10]{ZDG:95}.


\section{Preliminaries} \label{sec:prelim}

This section revises the Youla parametrization and the lifting technique, which are
required for technical developments in the paper. Although both subjects are
well-studied in the literature, both require some less documented twists.

\subsection{Youla parametrization with prespecified central controller}

Parametrizations of all stabilizing controllers for a given LTI plant, known as the
Youla, or Youla-Ku\v cera, parametrizations, is a classical result, well documented in
the literature, see \cite[Ch.~12]{ZDG:95} and the references therein. The idea also
extends to time-varying systems \cite[Sec.\:9.A]{FS:82}. These parametrizations are
conventionally expressed in terms of a linear-fractional transformation of a free
stable parameter (dubbed the ``$Q$-parameter'') by some given ``generator,'' which is
a function of a coprime factorization of the plant. When state-space realizations are
involved, the central controller, the one corresponding to $Q=0$, has commonly the
observer-based structure.

It is less common to construct a parametrization centered on some given, ``nominal,''
stabilizing controller, which is not necessarily observer based. This possibility was
explored in \cite[\S III-B]{Z:81} in the case when this nominal stabilizing controller
is stable itself. An insight into how to expand a given controller was provided in
\cite[pp.\,546--548]{BBN:90}, but constructive procedures and completeness were only
discussed for stable plants with the zero nominal controller and for observer-based
nominal controllers. I am not aware of other discussions of this subject in the
literature. Still, this kind of parametrization is required for developments in the
next section. Thus, although the result might not be entirely new, it is proved below.
\begin{lemma} \label{l:youla0}
 Let $P$ be an LTI plant having a strictly proper transfer function. Assume that it is
 internally stabilized by an LTI finite-dimensional controller $K_0$. Then all linear
 internally stabilizing controllers can be characterized as $K=\LFTl{J_0}{Q}$ for any
 stable and causal $Q$ and
 \begin{equation} \label{eq:allKfromK0}
  J_0\coloneq\mmatrix{K_0&\tilde M_0^{-1}\\M_0^{-1}&-M_0^{-1}P(I-K_0P)^{-1}\tilde M_0^{-1}},
 \end{equation}
 where $M_0$, $N_0$, $\tilde M_0$, and $\tilde N_0$ are coprime factors of $K_0$ over
 $R\Hinf$, such that $K_0=\tilde M_0^{-1}\tilde N_0=N_0M_0^{-1}$.
\end{lemma}
\begin{IEEEproof}
 Because $K_0$ is stabilizing, there must exist coprime factorizations of the plant,
 $P=\tilde M_P^{-1}\tilde N_P=N_PM_P^{-1}$, such that
 \[
  \mmatrix{\tilde M_0&-\tilde N_0\\-\tilde M_P&\tilde N_P}\mmatrix{M_P&N_0\\N_P&M_0}=I.
 \]
 Indeed, by \cite[Lem.~5.10]{ZDG:95} the stability of the closed-loop system implies
 that for any coprime factorizations of $P$, the systems $D\coloneq\tilde
 M_0M_P-\tilde N_0N_P$ and $\tilde D\coloneq\tilde M_PM_0-\tilde N_PN_0$ are
 bi-stable, i.e., such that $D,D^{-1},\tilde D,\tilde D^{-1}\in\Hinf$. Thus,
 $M_PD^{-1}$, $N_PD^{-1}$, $\tilde D^{-1}\tilde M_P$, and $\tilde D^{-1}\tilde N_P$
 are also coprime factors of $P$ and they do verify the equality above. It then
 follows from \cite[Sec.\:9.A]{FS:82} that all internally stabilizing controllers can
 be parametrized as $(N_0+M_PQ)(M_0+N_PQ)^{-1}$. The equivalence between this form and
 \eqref{eq:allKfromK0} follows by \cite[Lem.\:10.1]{ZDG:95} and the fact that
 $P(I-K_0P)^{-1}=N_P\tilde M_0$. Finally, as $P(\infty)=0$, the $(2,2)$ sub-block of
 $J_0(s)$ is strictly proper and the LFT in \eqref{eq:allKfromK0} is well posed for
 every causal $Q$ by \cite[Thm.\:4.1]{W:71}.
\end{IEEEproof}

By \cite[Lem.\,10.4(c)]{ZDG:95} the transformation $Q\mapsto K$ defined by
\eqref{eq:allKfromK0} is invertible, with $Q=\LFTu{J_0^{-1}}K$, where
\begin{equation} \label{eq:J0inv}
 J_0^{-1}=\mmatrix{P&M_0-PN_0\\\tilde M_0-\tilde N_0P&-\tilde N_0(M_0-PN_0)},
\end{equation}
and is well posed for any causal $K$, again by \cite[Thm.\:4.1]{W:71}.

\begin{remark}[connection with \cite{Z:81}]
 The parametrization of Lemma~\ref{l:youla0} can be rewritten as
 \[
  K=\LFTlL{\mmatrix{K_0&I\\I&-P(I-K_0P)^{-1}}}{\hat Q},
 \]
 where $\hat Q\coloneq\tilde M_0^{-1}QM_0^{-1}$. If $K_0$ is stable, both $\tilde M_0$
 and $M_0$ are bi-stable and can thus be absorbed into the $Q$-parameter. The
 parametrization then reduces to the case discussed in \cite{Z:81}. Yet unstable poles
 of $\tilde M_0^{-1}$ and $M_0^{-1}$, which are unstable poles of $K_0$, extend
 admissible $\hat Q$'s to a class of unstable systems. \Endoftheorem
\end{remark}

A state-space realization of the generator of all stabilizing controllers in
Lemma~\ref{l:youla0}, $J_0$, can also be derived. To this end, bring in stabilizable
and detectable realizations
\begin{equation} \label{eq:PandK0}
 P(s)=\mmatrix[c|c]{A&B_u\\\hline C_y&0}
 \quad\text{and}\quad
 K_0(s)=\mmatrix[c|c]{A_0&B_0\\\hline C_0&D_0}
\end{equation}
any pick any $F_0$ and $L_0$ such that $A_0+B_0F_0$ and $A_0+L_0C_0$ are Hurwitz.
Coprime factors of $K_0$ can then be constructed as in \cite[Thm.\:5.9]{ZDG:95}, which
eventually yields
\begin{align}
 J_0(s)
  &=\mmatrix[ccc|cc]{A_0&0&0&B_0&-L_0\\0&A_0&B_0C_y&0&-L_0\\0&B_uC_0&A+B_uD_0C_y&0&B_u\\\hline C_0&0&0&D_0&I\\-F_0&F_0&-C_y&I&0} \label{eq:J0ss} \\
  &=\mmatrix[c|cc]{A_0&B_0&-L_0\\\hline C_0&D_0&I\\-F_0&I&0}+\mmatrix{0&0\\0&J_\text a(s)} \nonumber
\end{align}
with stable
\[
 J_\text a(s)\coloneq\mmatrix[cc|c]{A_0&B_0C_y&-L_0\\B_uC_0&A+B_uD_0C_y&B_u\\\hline F_0&-C_y&0}.
\]

The state dimension of $J_0$ in \eqref{eq:J0ss} is in general higher than that of
$K_0$. For instance, consider the static feedback case, $K_0(s)=D_0$ for some $D_0$
such that the matrix $A+B_uD_0C_y$ is Hurwitz. Then
\begin{equation} \tag{\ref{eq:J0ss}$'$} \label{eq:J0ss:s}
 J_0(s)=\mmatrix[c|cc]{A+B_uD_0C_y&0&B_u\\\hline0&D_0&I\\-C_y&I&0},
\end{equation}
which is dynamic. In the observer-based case, where $K_0(s)=-F(sI-A-B_uF-LC_y)^{-1}L$
for some $F$ and $L$ such that $A+B_uF$ and $A+LC_y$ are Hurwitz, the dimension of
$J_0$ is not increased. It can be verified that the choices $F_0=C_y$ and $L_0=-B_u$
result then in $J_\text a=0$ and the parametrization with
\begin{equation} \tag{\ref{eq:J0ss}$''$} \label{eq:J0ss:ob}
 J_0(s)=\mmatrix[c|cc]{A+B_uF+LC_y&-L&B_u\\\hline F&0&I\\-C_y&I&0},
\end{equation}
as in \cite[Thm.~12.8]{ZDG:95}. For a general $K_0$, we may aim at picking admissible
$F_0$ and $L_0$ for which the order of $J_\text a$ is minimal.

\subsection{Lifting technique} \label{sec:lift}

The idea of lifting is to convert analog signals to discrete sequences of functions
operating over finite time intervals. Although mostly used to deal with systems with a
constant sampling rate, see \cite[Ch.\,10]{ChFr:95} and the references therein,
extensions of the technique to time-varying rates is effortless, at least at the level
required in this paper.


\begin{figure}[!tbp]
 \centering\includegraphics[scale=0.78]{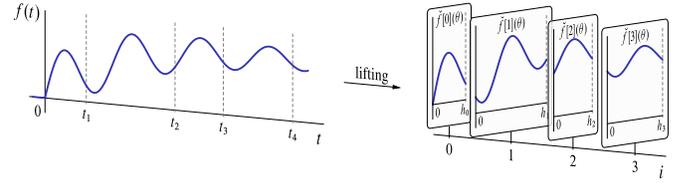}
 \caption{Lifting transformation with nonuniform time axis partition} \label{fig:lift}
\end{figure}

Consider a sequence of time instances $\{t_i\}_{i\in\mathbb Z^+}$ such that
$0=t_0<t_1<t_2<\cdots$. Then any analog signal $f:\mathbb R^+\to\mathbb R^n$ can be
equivalently cast as a sequence of functions $\{\vl f[i]\}_{i\in\mathbb Z^+}$ such
that $\vl f[i]:[0,h_i)\to\mathbb R^n$ is defined according to
\[
 \vl f[i](\theta)=f(t_i+\theta),\qquad i\in\mathbb Z^+,\theta\in[0,h_i)
\]
where $h_i\coloneq t_{i+1}-t_i$ is the length of the $i$th interval. The discrete
sequence $\{\vl f[i]\}$ is said to be the \emph{lifting} of the analog signal $f(t)$
with respect to the $t$-axis partition by $\{t_i\}$. See Fig.\,\ref{fig:lift} for a
visualization of this transformation.

\begin{figure}[!b]
 \centering
  \hspace*{\stretch1}
  \subfigure[In the time domain]{\label{fig:kerct}\includegraphics[scale=0.55]{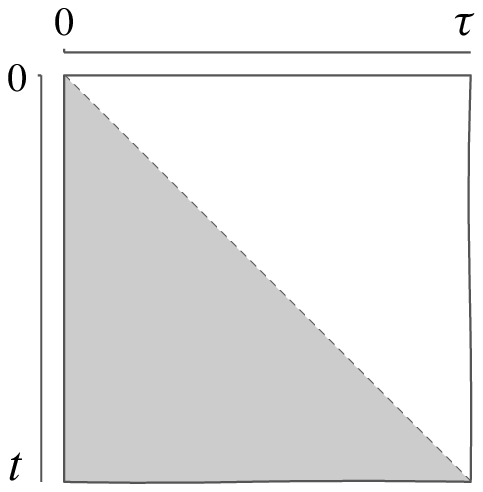}}
  \hspace*{\stretch2}
  \subfigure[In the lifted domain]{\label{fig:kerld}\includegraphics[scale=0.55]{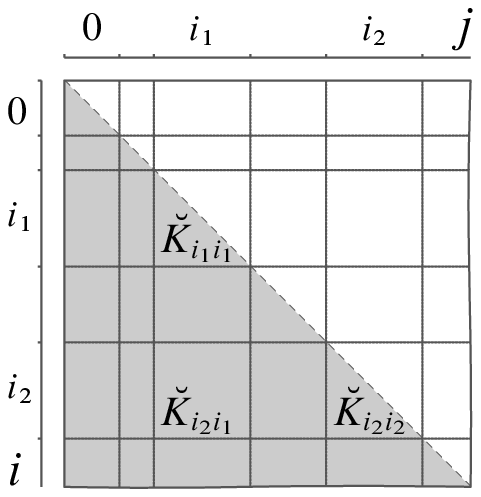}}
  \hspace*{\stretch2}
  \subfigure[The static part of $\vl K$]{\label{fig:kerlds}\includegraphics[scale=0.55]{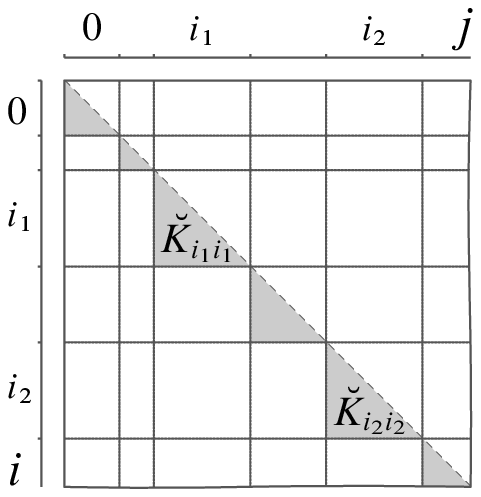}}
  \hspace*{\stretch1}
 \caption{Impulse responses of causal controllers} \label{fig:Kkers}
\end{figure}

Any continuous-time system can then be lifted by lifting its input and output signals,
resulting in a discrete-time system with infinite-dimensional input\,/\,output spaces.
To be specific, consider a causal controller $K:y\mapsto u$ described by the kernel
representation
\begin{equation} \label{eq:Gker}
 u(t)=\int_0^tk(t,\tau)y(\tau)\D\tau
\end{equation}
for an associated distribution $k(t,\tau)$ (impulse response) such that $k(t,\tau)=0$
whenever $t<\tau$. The impulse response may be visualizing as shown in
Fig.\,\ref{fig:kerct}, where the unshaded area represents zero values.
Relation \eqref{eq:Gker} can be rewritten in the lifted domain as
\[
 \vl u[i](\theta)
  =\sum_{j=0}^i\int_0^{h_j}\!k(t_i+\theta,t_j+\sigma)\vl y[j](\sigma)\D\sigma
  \eqcolon\biggl(\sum_{j=0}^i\vl K_{ij}\vl y[j]\biggr)(\theta)
\]
This relation describes a discrete linear system, denote it $\vl K$, whose kernel
(impulse response) $\vl K_{ij}$ at each $i,j$ is an integral operator mapping
functions on $[0,h_j)$ to functions on $[0,h_i)$. In terms of the kernel in
Fig.\,\ref{fig:kerct}, this transformation may be viewed as merely chopping the $t$-
and $\tau$-axes into pieces according to $\{t_i\}$. The result, shown in
Fig.\,\ref{fig:kerld}, can then be thought of as a form of system matrix as in
\cite[Sec.\:4.1]{ChFr:95}.
  
The ``diagonal'' elements $\vl K_{ii}$ of a lifted impulse response are called its
\emph{feedthrough parts}. At each $i\in\mathbb Z^+$ they are integral operators on
$[0,h_i)$ representing the direct connection between $\vl y[i]$ and $\vl u[i]$.  Given
a lifted system $\vl K$, by its \emph{static part} we understand the lifted system,
whose kernel is $\vl K_{ij}\delta_{ij}$, where $\delta_{ij}$ is the Kronecker delta.
At each $i$ this static part acts as $\vl u[i]=\vl K_{ii}\vl y[i]$, which corresponds
to the diagonal system matrix depicted in Fig.\,\ref{fig:kerlds}. The following
result, which is straightforward to verify, will be required in the sequel:
\begin{lemma} \label{l:spls}
 Let $G$ be an LTI system with the state-space realization $(A,B,C,D)$ and let $\vl G$
 be its lifting with respect to the time axis partition by $\{t_i\}$. Then the static
 part of $\vl G$ is the lifting of the continuous-time system $\zeta\mapsto\xi$
 verifying
 \[
  \begin{aligned}
   \dot x(t)&=Ax(t)+B\zeta(t),\quad x(t_i)=0 \\
   \xi(t)&=Cx(t)+D\zeta(t)
  \end{aligned}
 \]
 for all $t\in\mathbb R^+$ and $i\in\mathbb Z^+$.
\end{lemma}


\section{When is a Controller Sampled-Data\,?} \label{sec:sdc}

The redesign approach of this paper hinges on converting analog controllers to
sampled-data ones via \emph{constraining} the former. The first step to this end is to
understand how to characterize causal sampled-data controllers of the form in
Fig.\,\ref{fig:sdc} among linear operators mapping the measurement signal $y$ into the
control signal $u$.
  
\begin{figure}[!t]
 \centering
  \hspace*{\stretch1}
  \subfigure[In the time domain]{\label{fig:sdkerct}\includegraphics[scale=0.55]{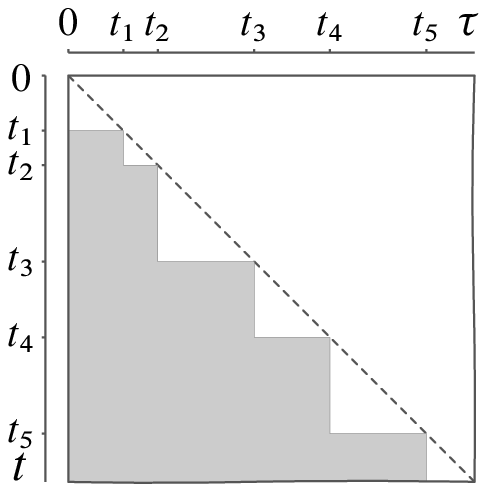}}
  \hspace*{\stretch2}
  \subfigure[In the lifted domain]{\label{fig:sdkerld}\includegraphics[scale=0.55]{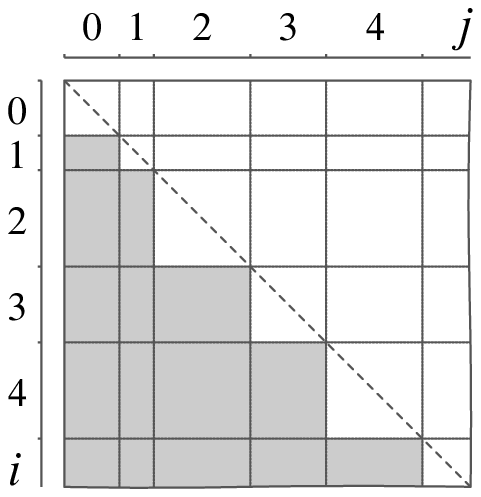}}
  \hspace*{\stretch1}
 \caption{Impulse responses of causal sampled-data controllers} \label{fig:sdKkers}
\end{figure}

An important role in the reasonings below is played by the fact that the sampler and
hold in Fig.\,\ref{fig:sdc} are not fixed (and neither are the dimensions of the
discrete signals $\vd y$ and $\vd u$). What should then be understood by a
sampled-data controller? The picture appears to be easier to grasp via causality of
the mapping $y\mapsto u$. Indeed, the very presence of the sampling operation inside
the controller should imply that between two subsequent sampling instances $u$ has no
new information about $y$, no matter what A/D and D/A converters are used. In other
words, $u(t)$ for all $t\in(t_i,t_{i+1})$ may be based on $y(\tau)$ for $\tau\le t_i$
only. \emph{Any controller satisfying this causality constraint will be regarded as an
  admissible one.} In terms of the kernel representation \eqref{eq:Gker},
admissibility then requires that
\begin{equation} \label{eq:sdkerc}
 \text{$k(t,\tau)=0$ whenever $\tau>\max_{t_i\le\,t}t_i$}.
\end{equation}
This yields a staircase, instead of triangle, constraint on the impulse response, as
shown in Fig.\,\ref{fig:sdkerct}.

Constraint \eqref{eq:sdkerc} might not be convenient to incorporate into design
procedures though, especially if the employed approach does not use the impulse
response directly. The constraint, however, is substantially simplified if translated
to the lifted domain associated with $\{t_i\}$. The stairs in Fig.\,\ref{fig:sdkerct}
fit then into the partition of the time axis, resulting in the system matrix in
Fig.\,\ref{fig:sdkerld}. This suggests that \eqref{eq:sdkerc} translates to the lifted
domain as \emph{strict causality}, i.e.\ the constraint that the feedthrough parts are
zero. The following result, which may be viewed as an extension of
\cite[Thm.\:1]{MR:97} to systems with non-uniform sampling and without the finite
dimensionality assumption about $K$, formalizes this observation:
\begin{theorem} \label{th:sdld}
 Let $\vl K$ be a linear system in the lifted domain with respect to the time axes
 partition by $\{t_i\}$. $\vl K$ is the lifting of a causal sampled-data system as in
 Fig.\,\ref{fig:sdc} with the sampling instances $\{t_i\}$ iff $\vl K$ is strictly
 causal, i.e.\ $\vl K_{ii}=0$, $\forall i\in\mathbb Z^+$.
\end{theorem}
\begin{IEEEproof}
 Follows by lifting \eqref{eq:sdkerc}.
\end{IEEEproof}

The strict causality is a more convenient system-theoretic notion to handle in various
controller design approaches than constraint \eqref{eq:sdkerc}. This is the reason to
introduce lifting.


\section{Stability-Preserving Redesign} \label{sec:stab}

Consider an LTI plant $P$. Without loss of generality, assume that its transfer
function $P(s)$ is strictly proper (this simplifies technicalities but can be easily
relaxed, see \cite[p.\,454]{ZDG:95}). Let a causal LTI controller $K_0$ internally
stabilize\footnote{The stability of a linear system $G$ is understood throughout the
  paper as its boundedness as an operator $L^2(\mathbb R^+)\to L^2(\mathbb R^+)$. In
  most cases the results remain unchanged if $L^2(\mathbb R^+)$ is replaced with
  $L^p(\mathbb R^+)$ for any $p\ge1$.} $P$ and $\{t_i\}_{i\in\mathbb Z^+}$ be a
sequence of time instances such that
\[
 0=t_0<t_1<\cdots<t_i<\cdots,\quad\text{with $\lim_{i\to\infty}t_i=\infty$}.
\]
The problem studied in this section is to approximate $K_0$ by a linear causal
sampled-data controller with the sampling instances $\{t_i\}$, so that the closed-loop
stability is preserved. By causal we understand a sampled-data controller as in
Fig.\,\ref{fig:sdc}, where $\mathcal S$ produces discrete signals $\vd y[i]$ at each
$t_i$ on the basis of measurements $y(t)$ for $t<t_i$, $\vd K$ is causal, and
$\mathcal H$ shapes the control signal $u(t)$ in $t\in[t_i,t_{i+1})$ on the basis of
discrete signals $\vd u[j]$ for $j\le i$. We assume hereafter that the sampling
instances $t_i$ are not known a priory, but the length of the intersample intervals
$h_i\coloneq t_{i+1}-t_i$ is uniformly bounded.

\subsection{Solution in the lifted domain} \label{sec:slds}

By Lemma~\ref{l:youla0}, $K_0$ generates the whole family of linear stabilizing
controllers, $K=\LFTl{J_0}Q$ for a given $J_0$, which is an augmentation of $K_0$, and
arbitrary stable and causal $Q$. Clearly, any stabilizing sampled-data controller must
belong to this family. It is therefore pertinent to understand, what conditions should
be imposed on $Q$ to produce sampled-data $\LFTl{J_0}Q$. The latter question, in turn,
is convenient to address in the lifted domain, where a handy characterization of
sampled-data controller exists, see Theorem~\ref{th:sdld}.

In the lifted domain, the controller parametrization reads $\vl K=\LFTl{\vl J_0}{\vl
  Q}$, where $\vl J_0$ and $\vl Q$ are the lifted versions of $J_0$ and $Q$,
respectively, with an arbitrary stable $\vl Q$ such that its feedthrough terms $\vl
Q_{ii}$ are causal. This LFT is then always well posed. Theorem~\ref{th:sdld} says
that $\vl K$ is the lifting of a sampled-data system iff its feedthrough terms $\vl
K_{ii}=0$ for all $i\in\mathbb Z^+$. The feedthrough terms of $\vl K$ depend only on
those of $\vl J_0$ and $\vl Q$ (because of their causality), i.e.\ $\vl
K_{ii}=\LFTl{\vl J_{0,ii}}{\vl Q_{ii}}$ for every $i$. Then, by
\cite[Lem.\,10.4(c)]{ZDG:95}, $\vl Q_{ii}=\LFTu{\vl J_{0,ii}^{-1}}{\vl K_{ii}}$.
Hence, for every $i$ we have that
\[
 \vl K_{ii}=0\iff\vl Q_{ii}=\vl Q_{0,ii}\coloneq\mmatrix{0&I}\vl J_{0,ii}^{-1}\mmatrix{0\\I}.
\]
This condition completely determines the feedthrough terms of $\vl Q$ and does not
affect the rest of it, which is handy.

Two straightforward, yet nevertheless important, observations are in order here.
First, $\vl Q_{0,ii}$ defined above is causal, because so is the continuous-time
system $J_0^{-1}$. Second, the static lifted system $\vl Q_\text{stat}$, whose impulse
response operators
\[
 \vl Q_{\text{stat},ij}=\ccases{\vl Q_{0,ii}&\text{if $j=i$}\\0&\text{otherwise}}
\]
is stable, as it is the lifting of an LTI system whose state resets at every $t_i$
with uniformly bounded\footnote{The uniform boundedness is actually required only if
  the $(2,2)$ sub-block of $J_0^{-1}$ in unstable. If this system is stable, which
  happens iff $P$ is itself stable (cf.\ \eqref{eq:J0inv}), the result holds for any
  $\{t_i\}$.} $t_{i+1}-t_i$. Consequently, any admissible $\vl Q$ can be presented as
$\vl Q=\vl Q_\text{stat}+\vl Q_\text{sd}$ for a strictly causal $\vl Q_\text{sd}$,
which is thus the lifting of a sampled-data system, and $\vl Q$ is stable iff $\vl
Q_\text{sd}$ is stable.

The discussion above can be summarized as follows:
\begin{lemma} \label{l:stabsdld}
 All causal stabilizing sampled-data controllers in the lifted domain can be
 parametrized as
 \[
  \vl K_\text{sd}=\LFTlL{\vl J_0}{\vl Q_\text{stat}+\vl Q_\text{sd}}
 \]
 for an arbitrary strictly causal stable $\vl Q_\text{sd}$, where $\vl Q_\text{stat}$
 is the static part of the $(2,2)$ sub-block of $\vl J_0^{-1}$.
\end{lemma}

\subsection{Solution in the continuous-time domain} \label{sec:stds}

Although treating the problem in the lifted domain is simple conceptually, it does not
result in a transparent solution. Our next step is thus to ``peel off'' the
lifted-domain result of Lemma~\ref{l:stabsdld}, i.e.\ to transform it back the time
domain, where the structure of the resulting controllers is clear.

To this end, let
\[
 J_0(s)=\mmatrix[c|cc]{A_J&B_{J1}&B_{J2}\\\hline C_{J1}&D_0&I\\C_{J2}&I&0},
\]
(concrete expressions of the parameters of this realization in terms of realizations
of $P$ and $K_0$ are given by \eqref{eq:J0ss}). The following theorem, which is a
sampled-data version of the Youla-Ku\v cera parametrization with unrestricted sampler
and hold, is then the main result of this section:
\begin{theorem} \label{th:Ksd}
 All causal stabilizing sampled-data controllers can be characterized as the
 interconnection of the sensor side ``pre-processor''
 \begin{subequations} \label{eq:Ksd}
 \begin{equation} \label{eq:Ksd:s}
  \dot x_\text s(t)=A_Jx_\text s(t)+B_{J1}y(t)+B_{J2}\bigl(u(t)-u_\text s(t)\bigr),
 \end{equation}
 where $u_\text s=C_{J1}x_\text s+D_0y$, and the ``post-processor''
 \begin{align}
  \dot x_\text a(t)&=(A_J-B_{J1}C_{J2})x_\text a(t)+B_{J2}\eta(t) \label{eq:Ksd:a} \\
  u(t)&=(C_{J1}-D_0C_{J2})x_\text a(t)+\eta(t) \label{eq:Ksd:u}
 \end{align}
 at the actuation side, connected via their sampled states as
 \begin{equation} \label{eq:Ksd:s2a}
  x_\text a(t_i)=x_\text s(t_i)
 \end{equation}
 and the signal $\eta=Q_\text{sd}(C_{J2}x_\text s+y)$, where $Q_\text{sd}$ is an
 arbitrary causal and stable sampled-data system.
 \end{subequations}
\end{theorem}
\begin{IEEEproof}
 The state-space realization of $J_0^{-1}$ is obtained by \cite[Lem.\:3.15]{ZDG:95}.
 Using Lemma~\ref{l:spls}, we then end up with $\vl
 Q_\text{stat}:\vl\epsilon\mapsto\vl\eta_Q$ as the lifting of
 \begin{equation} \label{eq:Qstatss}
  \begin{aligned}
   \dot x_Q(t)&=A_J^\times x_Q(t)-B_{J12}\epsilon(t),\quad x_Q(t_i)=0 \\
   \eta_Q(t)&=C_{J12}x_Q(t)-D_0\epsilon(t)
  \end{aligned}
 \end{equation}
 where $A_J^\times\coloneq A_J-B_{J1}C_{J2}-B_{J2}C_{J1}+B_{J2}D_0C_{J2}$,
 \begin{equation} \label{eq:BC12}
  B_{J12}\coloneq B_{J1}-B_{J2}D_0
  \quad\text{and}\quad
  C_{J12}\coloneq C_{J1}-D_0C_{J2}.
 \end{equation}
 Denoting by $\eta$ the output of $Q_\text{sd}$ and by $\epsilon$ the second output of
 $J_0$, the dynamics of $J_0$ read
 \[
  \begin{aligned}
   \dot x_J(t)&=A_Jx_J(t)+B_{J1}y(t)+B_{J2}\bigl(\eta(t)+\eta_Q(t)\bigr) \\
   u(t)&=C_{J1}x_J(t)+D_0y(t)+\eta(t)+\eta_Q(t) \\
   \epsilon(t)&=C_{J2}x_J(t)+y(t)
  \end{aligned}
 \]
 (the second input of $J_0$ is the sum of the outputs of $Q_\text{stat}$ and
 $Q_\text{sd}$). Combining this realization with \eqref{eq:Qstatss}, eliminating
 $\eta_Q$, and carrying out a state transformation yields \eqref{eq:Ksd} with $x_\text
 s=x_J$ and $x_\text a=x_J+x_Q$.
\end{IEEEproof}

The signal $u_\text s$ in pre-processor \eqref{eq:Ksd:s} may be thought of as an
emulation of the output of the analog controller $K_0$, which equals $C_{J1}x_J+D_0y$.
The pre-processor resembles then the state observer for $J_0$.  The only difference is
that the calculated output, $u_\text s$, is now compared with the actual control
signal, $u$, produced by another system, via the sampling operation
\eqref{eq:Ksd:s2a}.

The central controller, the one with $Q_\text{sd}=0$ (and $\eta=0$), can be presented
in the form shown in Fig.\,\ref{fig:sdc}. To describe its components, introduce the
matrix functions
\begin{multline*}
 \mmatrix{\Lambda_{11}(\theta)&\Lambda_{12}(\theta)\\0&\Lambda_{22}(\theta)}\\\coloneq\exp\biggl(\mmatrix{A_J-B_{J2}C_{J1}&B_{J2}C_{J12}\\0&A_J-B_{J1}C_{J2}}\theta\biggr)
\end{multline*}
with $\Lambda_{11}(\theta)=\E^{(A_J-B_{J2}C_{J1})\theta}$,
$\Lambda_{22}(\theta)=\E^{(A_J-B_{J1}C_{J2})\theta}$, and
\begin{equation} \label{eq:Lam12}
 \Lambda_{12}(\theta)=\int_0^\theta\Lambda_{11}(\theta-\sigma)B_{J2}C_{J12}\Lambda_{22}(\sigma)\D\sigma
\end{equation}
(by Van Loan's formulae, see e.g.\ \cite[Lem.\:10.5.1]{ChFr:95}). Then:
\begin{corollary} \label{cor:Ksd}
 The ``central'' controller of Theorem~\ref{th:Ksd} can be implemented as the
 sampled-data controller in Fig.\,\ref{fig:sdc} with the generalized sampler (A/D
 converter) $\mathcal S:y\mapsto\vd y$
 \begin{subequations} \label{eq:K0sd}
 \begin{equation} \label{eq:S0}
  \vd y[i+1]=\int_0^{h_i}\E^{(A_J-B_{J2}C_{J1})(h_i-\sigma)}B_{J12}y(t_i+\sigma)\D\sigma,
 \end{equation}
 the discrete-time controller $\vd K:\vd y\mapsto\vd u$
 \begin{equation} \label{eq:K0d}
  \vd u[i+1]=(\Lambda_{11}(h_i)+\Lambda_{12}(h_i))\vd u[i]+\vd y[i+1],
 \end{equation}
 and the generalized hold (D/A converter) $\mathcal H:\vd u\mapsto u$
 \begin{equation} \label{eq:H0}
  u(t_i+\theta)=C_{J12}\,\E^{(A_J-B_{J1}C_{J2})\theta}\vd u[i],
 \end{equation}
 \end{subequations}
 where $B_{J12}$ and $C_{J12}$ are defined by \eqref{eq:BC12}.
\end{corollary}
\begin{IEEEproof}
 Rewrite \eqref{eq:Ksd:s} as
 \[
  \dot x_\text s(t)=(A_J-B_{J2}C_{J1})x_\text s(t)+B_{J12}y(t)+B_{J2}u(t),
 \]
 so that
 \begin{multline*}
  x_\text s(t_{i+1})=\Lambda_{11}(h_i)x_\text s(t_i)+\smash[t]{\int_0^{h_i}}\Lambda_{11}(h_i-\sigma)\\\times\bigl(B_{J12}y(t_i+\sigma)+B_{J2}u(t_i+\sigma)\bigr)\D\sigma.
 \end{multline*}
 Now, \eqref{eq:S0}, the fact that $u(t_i+\sigma)=C_{J12}\Lambda_{22}(\sigma)x_\text
 s(t_i)$, which follows from \eqref{eq:Ksd:a}--\eqref{eq:Ksd:s2a} with $\eta=0$, and
 \eqref{eq:Lam12} yield that
 \[
  x_\text s(t_{i+1})=(\Lambda_{11}(h_i)+\Lambda_{12}(h_i))x_\text s(t_i)+\vd y[i+1].
 \]
 The result follows by introducing $\bar u[i]\coloneq x_\text s(t_i)$.
\end{IEEEproof}

Controller \eqref{eq:K0sd} is well suited to networked implementation. Sampler
\eqref{eq:S0} requires uninterrupted access to the measured output $y$ and should be
implemented on the sensor side. Hold \eqref{eq:H0} generates a complex waveform analog
control signal $u$, so it should be implemented on the actuator side. The exchange of
information between these parts, done via \eqref{eq:K0d}, may be intermittent. It can
be carried out either opportunistically, when network resources are available, or when
menacing deviations from predicted behavior are detected. In any case, the nominal
closed-loop system remains stable for any uniformly bounded sequence of sampling
intervals $\{h_i\}$.

The control signal $u$ generated by \eqref{eq:Ksd} is typically discontinuous because
of jumps in $x_\text a$ at $t=t_i$, cf.\ \eqref{eq:Ksd:s2a}. A workaround is to
parametrize the set of analog stabilizing controllers in the form
$K=F_\text{lp}\LFTl{\tilde J_0}Q$ for some low-pass $F_\text{lp}$. This can be done by
factoring $K_0=F_\text{lp}\tilde K_0$ and then applying Lemma~\ref{l:youla0} to
$\tilde K_0$ and the augmented plant $\tilde P=PF_\text{lp}$. In this case only
$\LFTl{\tilde J_0}Q$ is redesigned, so that the actual control signal is a filtered
version of \eqref{eq:Ksd:u}. This factorization is sometimes a part of the design
method, see Section~\ref{sec:ex} for an example.

\subsection{Special cases} \label{sec:specc}

To illustrate the structure of the controller derived above, consider in this
subsection some special cases. It is assumed throughout that the plant is given in
terms of its state-space realization \eqref{eq:PandK0}.

\subsubsection{Static $K_0$}

Let $K_0(s)=D_0$ for a $D_0$ such that the matrix $A+B_uD_0C_y$ is Hurwitz. Then
$J_0(s)$ is given by \eqref{eq:J0ss:s} and \eqref{eq:Ksd} can be rewritten as
\begin{subequations} \label{eq:Ksds}
\begin{align}
 \dot x_\text s(t)&=Ax_\text s(t)+B_uu(t)-B_uD_0\bigl(y(t)-C_yx_\text s(t)\bigr) \label{eq:Ksds:s} \\
 \dot x_\text a(t)&=Ax_\text a(t)+B_uu(t),\quad x_\text a(t_i)=x_\text s(t_i) \label{eq:Ksds:a} \\
 u(t)&=D_0C_yx_\text a(t)+\eta(t) \label{eq:Ksds:u}
\end{align}
\end{subequations}
The sensor-side part, \eqref{eq:Ksds:s}, is the standard full-order observer of the
plant state with the gain $L=B_uD_0$. The actuator-side part,
\eqref{eq:Ksds:a}--\eqref{eq:Ksds:u}, mimics then the dynamics of the closed-loop
system under the analog control law $u=D_0y+\eta$.

\subsubsection{Observer-based $K_0$}

In this case the generator of all stabilizing controllers, $J_0$, is given by
\eqref{eq:J0ss:ob}. Hence, \eqref{eq:Ksd:s} reads%
\begin{subequations} \label{eq:Ksdob}
\begin{align}
 \dot x_\text s(t)&=Ax_\text s(t)+B_uu(t)-L\bigl(y(t)-C_yx_\text s(t)\bigr), \label{eq:Ksdob:s} \\[-.75ex]
\intertext{\smash{which is again an observer, and \eqref{eq:Ksd:a}--\eqref{eq:Ksd:s2a} read}}\nonumber \\[-4.25ex]
 \dot x_\text a(t)&=Ax_\text a(t)+B_uu(t),\quad x_\text a(t_i)=x_\text s(t_i) \label{eq:Ksdob:a} \\
 u(t)&=Fx_\text a(t)+\eta(t). \label{eq:Ksdob:u}
\end{align}
\end{subequations}
In the intermittent sampling case, this controller structure was proposed in
\cite{GW:10}, although with no stability proof. Apparently, the first proof of the
closed-loop stability under this scheme was offered in \cite{LL:11}. In the constant
$h_i$ case, earlier proofs exist. If presented in form \eqref{eq:K0sd}, this is
exactly the optimal controller configuration of \cite[Thm.\:5.1]{MRP:99I}. The even
earlier result of \cite[Thm.\:3.1]{Tad:92} is also essentially the same system, sans
the absorption of $Q_\text{stat}$ into $J_0$. See also \cite[Ch.\:3]{GAM:14} for an
analysis of the same controller under the constant sampling rate and parametric plant
uncertainty.

Curiously, the redesigned static controller \eqref{eq:Ksds} is a special case of the
redesigned observer-based controller \eqref{eq:Ksdob}, under $L=B_uD_0$ and
$F=D_uC_y$. Consequently, the use of static controllers offers no advantage over
observer-based controllers in terms of simplicity for the proposed redesign procedure.

\subsection{Complexity reduction via $Q_\text{sd}$} \label{sec:complr}

The freedom in the choice of $Q_\text{sd}$ can be used to reduce the complexity of
the controller of Theorem~\ref{th:Ksd}. Consider, for example, the following
$Q_\text{sd}:y-C_{J2}x_\text s\mapsto\eta$:
\begin{equation} \label{eq:Qsd1}
 \begin{aligned}
  \dot x_\eta(t)&=A_\eta x_\eta(t),&x_\eta(t_i)=B_\eta(y(t_i)-C_{J2}x_\text s(t_i)) \\
  \eta(t)&=C_\eta x_\eta(t)
 \end{aligned}
\end{equation}
which is the cascade of the ideal sampler and a generalized hold as in \eqref{eq:H0},
just with different parameters. System \eqref{eq:Qsd1} is stable for any $A_\eta$,
$B_\eta$, and $C_\eta$, because it resets at every $t_i$. With this choice, the
actuation-side dynamics \eqref{eq:Ksd:a}--\eqref{eq:Ksd:u} read
\[
 \begin{aligned}
  \mmatrix{\dot x_\text a(t)\\\dot x_\eta(t)}&=\mmatrix{A_J-B_{J1}C_{J2}&B_{J2}C_\eta\\0&A_\eta}\mmatrix{x_\text a(t)\\x_\eta(t)} \\
  u(t)&=\mmatrix{C_{J1}-D_0C_{J2}&C_\eta}\mmatrix{x_\text a(t)\\x_\eta(t)}
 \end{aligned}
\]
with the following effect of \eqref{eq:Ksd:s} on them:
\[
 \mmatrix{x_\text a(t_i)\\x_\eta(t_i)}=\mmatrix{I\\-B_\eta C_{J2}}x_\text s(t_i)+\mmatrix{0\\B_\eta}y(t_i).
\]
If $C_\eta=C_{J1}-D_0C_{J2}$, then $u$ depends only on $\tilde x_\text a\coloneq
x_\text a+x_\eta$. If then $A_\eta=A_J^\times$ defined after \eqref{eq:Qstatss}, the
signal $\tilde x_\text a$ becomes independent of $x_\eta$ (can be seen by a similarity
transformation). As a result, we end up with essentially unchanged actuator-end
equations (just with $\eta=0$) and with the new interconnection
\begin{equation} \tag{\ref{eq:Ksd:s2a}$'$}\label{eq:Ksdr:s2a}
 x_\text a(t_i)=(I-B_\eta C_{J2})x_\text s(t_i)+B_\eta y(t_i).
\end{equation}
in place of \eqref{eq:Ksd:s2a}. We may then seek for $B_\eta$ that renders some modes
of \eqref{eq:Ksd:s}, which are the eigenvalues of $A_J-B_{J2}C_{J1}$, unobservable
through $I-B_\eta C_{J2}$. Unobservable dynamics may then be safely canceled, reducing
the order of \eqref{eq:Ksd:s}.

A possible procedure for carrying out such a reduction is as follows. Assume w.l.o.g.\ 
that $C_{J2}$ has full row rank. Let $V_2$ be a matrix such that $\im V_2$ is
$(A_J-B_{J2}C_{J1})$-invariant and $C_{J2}V_2$ is left invertible. Pick $B_\eta$ as
any solution of $B_\eta C_{J2}V_2=V_2$. In this case $\im V_2=\ker(I-B_\eta C_{J2})$,
which implies that $\im V_2$ is the unobservable subspace of the $(I-B_\eta
C_{J2},A_J-B_{J2}C_{J1})$. Hence, all modes of $A_J-B_{J2}C_{J1}|\im V_2$ are
unobservable through $I-B_\eta C_{J2}$ and can thus be canceled. The maximal reduction
is attained if there is an admissible $V_2$ such that $C_{J2}V_2$ is square.
 
The choice of $B_\eta$ is particularly simple in the static state-feedback case, which
corresponds to \eqref{eq:Ksds} with $C_y=I$ and $D_0=F$ for some $F$ such that
$A+B_uF$ is Hurwitz. With the choice $B_\eta=I$, equation \eqref{eq:Ksdr:s2a} reads
$\tilde x_\text a(t_i)=x(t_i)$, which renders observer \eqref{eq:Ksds:s} redundant.
This yields the control law
\[
 u(t)=F\E^{(A+B_uF)(t-t_i)}x(t_i),\quad\forall t\in[t_i,t_{i+1})
\]
which effectively reproduces the algorithm of \cite{MA:04} (see also
\cite[Ch.\:5]{GAM:14}) and \cite{LL:10} (the latter also adds the effect of a
piece-wise constant disturbance estimate to the generated $u$).


\section{Performance-Guaranteeing Redesign} \label{sec:perfn}

The procedure of Section~\ref{sec:stab} produces a family of stabilizing sampled-data
controllers from a given analog controller $K_0$. Of this family one would naturally
prefer a controller that is close to $K_0$, in whatever sense. This section studies
situations when the closeness between $K_0$ and its sampled-data approximation is
measured in terms of the attained closed-loop performance.

\begin{figure}[tb]
 \centering
  \hspace*{\stretch1}
  \subfigure[Analog controller]{\label{fig:spct}\includegraphics[scale=0.667]{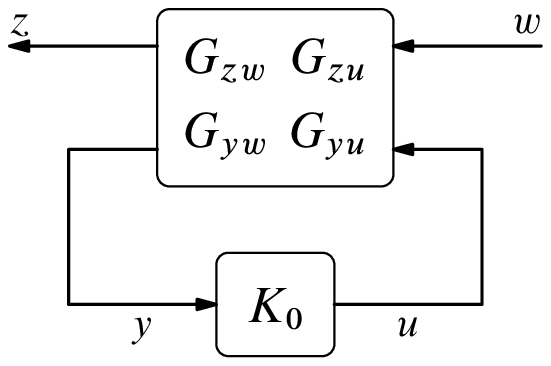}}
  \hspace*{\stretch2}
  \subfigure[Sampled-data controller]{\label{fig:spsd}\includegraphics[scale=0.667]{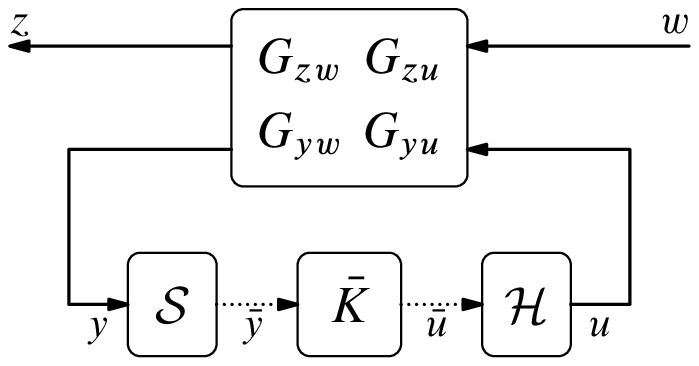}}
  \hspace*{\stretch1}
 \caption{Standard problems} \label{fig:sp}
\end{figure}

To this end, the setup is extended to the so-called ``standard problem'' of the form
depicted in Fig.\,\ref{fig:spct}. The performance of this system is quantified by a
norm, either \Htwo\ or \Hinf, of the closed-loop system $T_{zw}\coloneq\LFTl G{K_0}$
from $w$ to $z$. It is assumed that $K_0$ guarantees certain performance level and the
goal is to find a sampled-data controller that can deliver a comparable performance
level for the setup in Fig.\,\ref{fig:spsd}.

\begin{remark}[viewpoint]
 The problems addressed in this section might also be viewed as merely the design of
 (sub)\,optimal sampled-data controllers for intermittent sampling. But optimality
 might make little engineering sense per se. Rather, it is a powerful tool to design
 ``good'' analog controllers. For that reason, solving the very same optimization
 problem for a sampled-data controller is treated here as a tool of redesigning a
 chosen analog controller $K_0$. \Endoftheorem
\end{remark}

Throughout this section, we assume that
\[
 G(s)=\mmatrix{G_{zw}(s)&G_{zu}(s)\\G_{yw}(s)&G_{yu}(s)}=\mmatrix[c|cc]{A&B_w&B_u\\\hline C_z&0&D_{zu}\\C_y&D_{yw}&0}
\]
and that the standard assumptions \cite[p.\:384]{ZDG:95} are satisfied (including the
normalizations $D_{zu}'D_{zu}=I$ and $D_{yw}D_{yw}'=I$). The solution procedure is
again to start with a parametrization, now of all suboptimal analog controllers, and
then seek for a ``least harmful'' $Q$-parameter for which the resulting controller is
a sampled-data one.

\subsection{\Htwo\ performance} \label{sec:H2}

Let $K_0$ be the \Htwo-optimal controller for the problem in Fig.\,\ref{fig:spct} and
$\{t_i\}$ be a sequence of sampling instances. The problem studied below is to find
the optimal sampled-data controller, of the form depicted in Fig.\,\ref{fig:sdc}, for
the same generalized plant.

The \Htwo\ norm of a linear system can be roughly viewed as the $L^2(\mathbb
R^+)$-norm of its impulse response. In the LTI case, it is sufficient to consider the
response to the impulse applied at $t=0$, which leads to the conventional definition
\cite[p.\:98]{ZDG:95}. The response of time-varying systems to impulses applied at
different time instances might differ. A way to generalize the notion of the \Htwo\ 
norm to such systems is via averaging. Namely, let $G$ be a linear system described by
\eqref{eq:Gker}. Then we may define (see e.g.\ \cite{IW:92scl} or
\cite[\S2.1.2]{IK:01}) its \Htwo\ norm as
\begin{equation} \label{eq:H2norm}
 \norm{G}_2^2\coloneq\lim_{T\to\infty}\frac1T\int_0^T\int_\tau^\infty\normf{g(t,\tau)}^2\,\D t\D\tau,
\end{equation}
where $\normf\cdot$ denotes the Frobenius matrix norm. This quantity may also be
thought of as the average output variance if the input is a zero mean white noise
process. In general, \eqref{eq:H2norm} is a semi-norm, although in some special cases,
like periodic systems, it is a norm. It reduces to the standard definition if $G$ is
time invariant.

The main result of this sub-section is formulated below:
\begin{theorem} \label{th:H2}
 Let the analog \Htwo\ problem associated with the system in Fig.\,\ref{fig:spct} be
 well posed and $F$ and $L$ be the state-feedback and filter gains associated with
 this problem. Then the optimal \Htwo\ performance attainable by sampled-data
 controllers for a given sequence of sampling instances $\{t_i\}$ is
 \[
  \gamma_{\{t_i\}}^2=\gamma_0^2+\lim_{i\to\infty}\frac1{t_i}\sum_{j=0}^{i-1}\int_0^{h_j}\!\int_0^{h_j-\tau}\normf{F\E^{At}L}^2\,\D t\D\tau,
 \]
 where $\gamma_0$ is the optimal \Htwo\ performance attainable by analog controllers.
 The \Htwo\ performance attained by the sampled-data controller given by
 \eqref{eq:Ksdob} with $\eta=Q_\text{sd}(y-C_yx_\text s)$ is then
 $\norm{T_{zw}}_2^2=\gamma_{\{t_i\}}^2+\norm{Q_\text{sd}}_2^2$.
\end{theorem}
\begin{IEEEproof}
 See Appendix.
\end{IEEEproof}

Note that the optimal sampled-data controller is not unique. Because \eqref{eq:H2norm}
is a semi-norm, there are nonzero $Q_\text{sd}$ such that $\norm{Q_\text{sd}}_2=0$.
Any such $Q_\text{sd}$ produces an optimal controller.

An intriguing question is under what sampling pattern $\{t_i\}$ the attainable
performance is minimal. Of course, this question makes sense only if the ``average''
sampling period is fixed. Another assumption that should be made in this respect is
that the sampling pattern is periodic. Otherwise, an alternation of any finite number
of sampling instances $t_i$ has no effect on $\gamma_{\{t_i\}}$. Thus, assume that
there is an $N$ such that $h_{i+N}=h_i$ for all $i\in\mathbb Z^+$ and that
\begin{equation} \label{eq:averh}
 h_\text{av}\coloneq\frac1N\sum_{i=0}^{N-1}h_i=\frac{t_N}N
\end{equation}
is fixed. In this case
\[
 \gamma_{\{t_i\}}^2=\gamma_0^2+\frac1{Nh_\text{av}}\sum_{j=0}^{N-1}\int_0^{h_j}\!\int_0^{h_j-\tau}\normf{F\E^{At}L}^2\,\D t\D\tau
\]
(and, as a matter of fact, the optimal $Q_\text{sd}=0$ is unique now). The optimal
sampling pattern is then given as follows:
\begin{proposition} \label{pr:opttkH2}
 If $K_0\ne0$, the unique optimal sampling pattern for a fixed $h_\text{av}$ in
 \eqref{eq:averh} and any $N\in\mathbb Z^+\setminus\{0\}$ is the uniform sampling,
 i.e.\ $h_i=h_\text{av}$ for all $i\in\mathbb Z^+$.
\end{proposition}
\begin{IEEEproof}
 First, $K_0(s)=-F(sI-A-B_uF-LC_y)^{-1}L=0$ iff $F(sI-A)^{-1}L=0$, which is readily
 verified via the Kalman canonical decomposition \cite[Thm.\:3.10]{ZDG:95}. Hence, the
 condition of the proposition guarantees that $F\E^{At}L\not\equiv0$ in any finite
 interval of $\mathbb R^+$.
 
 Let us start with the case of $N=2$. Sampling periods can then be parametrized as
 $h_0=h-\delta$ and $h_1=h+\delta$ for $\delta\in[-h,h]$ and the optimal performance
 is
 \[
  \gamma_{\{t_i\}}^2=\gamma_0^2+\frac{\gamma_1(h+\delta)+\gamma_1(h-\delta)}{2h},
 \]
 where
 \[
  \gamma_1(h)\coloneq\int_0^h\int_0^{h-\tau}\!\normf{F\E^{At}L}^2\,\D t\D\tau.
 \]
 It can be verified, using the Leibniz integral rule, that
 \[
  \frac{\D\gamma_1(h+\delta)}{\D\delta}=\int_0^{h+\delta}\normf{F\E^{At}L}^2\,\D t,
 \]
 so that
 \[
  \frac{\D\gamma_{\{t_i\}}^2}{\D\delta}=\frac1{2h}\int_{h-\delta}^{h+\delta}\normf{F\E^{At}L}^2\,\D t
 \]
 has the same sign as $\delta$ and is zero iff $\delta=0$. This proves the statement
 of the Proposition.
 
 Now consider the case of $N>2$. If not all $h_i$ are equal, we can always find a
 $j>1$ such that $h_{j-1}\ne h_j$. The replacement of $t_j$ with $(t_{j+1}+t_{j-1})/2$
 then decreases $\gamma_1(h_{j-1})+\gamma_1(h_j)$ and affects no other
 $\gamma_1(h_i)$. Hence, there always a pattern yielding a better performance. This
 procedure fails to reduce $\gamma_{\{t_i\}}$ only if all $h_i=h_\text{av}$, which
 completes the proof.
\end{IEEEproof}

Proposition~\ref{pr:opttkH2}, which establishes that the uniform sampling is
advantageous, appears to disagree with some earlier results. This aspect is clarified
in the following two remarks.
\begin{remark}[alternative choices of the \Htwo\ norm]
 A variable sampling rate scheme to improve the LQR performance in sampled-data
 systems was proposed in \cite{BB:14}. It is based on the rate of change of the
 optimal analog control signal and is optimal for 1-order systems. The problem studied
 in \cite{BB:14} is different from that studied here though. First, it assumes the
 zero-order hold and the ideal sampler. This is different, and more restrictive, from
 the setup with free hold and sampler. Second, and most importantly, the performance
 measure considered in \cite{BB:14} is different. The LQR optimization effectively
 minimizes the energy of the response to the impulse applied at $t=0$ only. In other
 words, it does not involve averaging. As follows from the proof of
 Theorem~\ref{th:H2}, if this philosophy were used in the \Htwo\ design for the system
 in Fig.\,\ref{fig:spct}, the optimal performance would be
 \[
  \norm{T_{zw}}_2^2=\gamma_0^2+\int_0^{h_0}\normf{F\E^{At}L}^2\,\D t.
 \]
 The obvious choice is then $t_1\to0$, which recovers the analog performance
 irrespective of the other sampling instances. But this design would make no practical
 sense. Another possibility, something between \eqref{eq:H2norm} and LQR, would be to
 consider
 \[
  \norm G_2^2\coloneq\lim_{i\to\infty}\frac1i\sum_{j=0}^{i-1}\int_{t_j}^\infty\normf{g(t,t_j)}^2\,\D t.
 \]
 Consider what happens with this choice when the sampling pattern is $2$-periodic. In
 that case,
 \[
  \norm{T_{zw}}_2^2=\gamma_0^2+\frac12\biggl(\int_0^{h+\delta}\normf{F\E^{At}L}^2\,\D t+\int_0^{h-\delta}\normf{F\E^{At}L}^2\,\D t\biggr)
 \]
 so that\vspace*{-2.5pt}
 \[
  \frac{\D\norm{T_{zw}}_2^2}{\D\delta}=\frac{\normf{F\E^{A(h+\delta)}L}^2-\normf{F\E^{A(h-\delta)}L}^2}2.
 \]
 Similarly to the proof of Proposition~\ref{pr:opttkH2}, this function equals zero at
 $\delta=0$. But this might neither be the only such point nor the point of the local
 minimum, depending on the parameters. For example, assume that the system is 1-order,
 i.e.\ $A$, $F$, and $L$ are scalars. In this case, the sign of the derivative of the
 optimal performance equals $\sign(\E^{A\delta}-\E^{-A\delta})$. Thus, if the system
 is unstable ($A>0$), the uniform sampling is still the best option. But if the system
 is stable ($A<0$), the uniform sampling is the worst scenario and the best option is
 to alternate short and long sampling intervals. If $A=0$, the sampling pattern is
 irrelevant. If $G$ has higher order dynamics, the optimal sampling pattern might be
 more complicated. \Endoftheorem
\end{remark}

\begin{remark}[realization vs.\ process]
 Another way to assign the sampling pattern is to use event-based mechanisms
 \cite{A:08,HJT:12}. Some results of this kind analyze the \Htwo\ performance. For
 example, the Lebesgue sampling strategy of \cite{AB:03} (see also
 \cite[Sec.\:3]{A:08}) may result in a significant relaxation of the average sampling
 rate (by a factor of $3$ in the case where $A=D_\bullet=0$ and
 $B_\bullet=C_\bullet=1$). The cause of this improvement may lie in the ability of
 event-based sampling to make use of the information about the effect of a
 \emph{particular realization} of $w$ on the system, rather than treating $w$ as a
 \emph{random process}. It may be interesting in this respect to investigate the
 possibility to use the signal $Q_\text{stat}(y-C_yx_\text s)$, with $Q_\text{stat}$
 as in \eqref{eq:QstatH2}, as the basis for event generation. This would be
 qualitatively different from existing event generation mechanisms as it involves
 low-pass filtering of the estimation error. This element may be useful in avoiding
 Zeno behavior \cite{HJT:12} and may lead to performance-justified switching, see the
 example in \S\ref{sec:pend}. \Endoftheorem
\end{remark}

\subsection{\Hinf\ performance} \label{sec:Hinf}

Unlike the \Htwo\ case, the \Hinf\ performance measure admits a clean and unambiguous
generalization to time-varying systems, as the $L^2(\mathbb R^+)$ induced norm. Denote
by $\gamma_\text{opt}\ge0$ the optimal \Hinf\ performance attainable for the standard
problem associated with Fig.\,\ref{fig:spct} by an analog controller. Let $K_0$ be the
central $\gamma$-suboptimal controller for a $\gamma>\gamma_\text{opt}$. This $K_0$
generates the whole family of $\gamma$-suboptimal controllers. The question asked
below is under what conditions on the sequence of sampling instances $\{t_i\}$ this
family contains a sampled-data controller of the form depicted in Fig.\,\ref{fig:sdc}.

To formulate the result, we need the Riccati equations
\begin{gather*}
 XA+A'X+C'_zC_z+\gamma^{-2}XB_wB'_wX-F'F=0, \\
 AY+YA'+B_wB'_w+\gamma^{-2}YC'_zC_zY-LL'=0,
\end{gather*}
where $F\coloneq-B'_uX-D'_{zu}C_z$ and $L\coloneq-YC'_y-B_wD'_{yw}$. The solutions $X$
and $Y$ are called stabilizing if the matrices $A_F\coloneq
A+\gamma^{-2}B_wB_w'X+B_uF$ and $A_L\coloneq A+\gamma^{-2}YC_z'C_z+LC_y$ are Hurwitz.
It is known \cite[Thm.\:16.4]{ZDG:95} that $\gamma>\gamma_\text{opt}$ iff the
stabilizing solutions exist and are such that $X\ge0$, $Y\ge0$, and
$\rho(YX)<\gamma^2$. We then have:
\begin{theorem} \label{th:Hinf}
 Let $\gamma>\gamma_\text{opt}$. Then there is a $\gamma$-suboptimal sampled-data
 controller for a given sequence of sampling instances $\{t_i\}$ iff there exists a
 solution to the differential Riccati equation
 \begin{multline*}
  \dot P(t)=AP(t)+P(t)A'\\+B_wB_w'+\gamma^{-2}P(t)C_z'C_zP(t),\quad P(0)=Y
 \end{multline*}
 such that $\rho(P(t)X)<\gamma^2$, $\forall t\in[0,h_i]$ and every $i\in\mathbb Z^+$.
 If the condition holds, a $\gamma$-suboptimal sampled-data controller is%
 \begin{subequations} \label{eq:KsdHinf}
 \begin{align}
  \dot x_\text s(t)&=A_Lx_\text s(t)-Ly(t)+(B_u+\gamma^{-2}YC_z'D_{zu})u(t), \label{eq:KsdHinf:s} \\
  \dot x_\text a(t)&=A_Fx_\text a(t),\quad x_\text a(t_i)=(I-\gamma^{-2}YX)^{-1}x_\text s(t_i) \label{eq:KsdHinf:a} \\
  u(t)&=Fx_\text a(t). \label{eq:KsdHinf:u}
 \end{align}
 \end{subequations}
\end{theorem}
\begin{IEEEproof}
 See Appendix.
\end{IEEEproof}

\begin{remark}[closed-loop stability] \label{rem:rob}
 The stability of the closed-loop system under the control law \eqref{eq:KsdHinf} is
 guaranteed only if the condition of Theorem~\ref{th:Hinf} holds for all $h_i$. This
 is in contrast to the \Htwo\ case, where the controller is stabilizing even if it
 does not guarantee a required performance level. \Endoftheorem
\end{remark}
\begin{remark}[generating disturbances]
 In terms of $\tilde x_\text s\coloneq(I-\gamma^{-2}YX)^{-1}x_\text s$ the sensor-side
 dynamics in \eqref{eq:KsdHinf:s} read
 \begin{multline*}
  \dot{\tilde x}_\text s(t)=A\tilde x_\text s(t)+B_w\tilde w_\gamma(t)+B_uu(t)\\-\tilde L\bigl(y(t)-C_y\tilde x_\text s(t)-D_{yw}\tilde w_\gamma(t)\bigr),
 \end{multline*}
 where $\tilde L\coloneq(I-\gamma^{-2}YX)^{-1}L$ and $\tilde
 w_\gamma\coloneq\gamma^{-2}B_w'X\tilde x_\text s$. This is the \Hinf\ estimator for
 the analog control signal $u=Fx$ in the presence of the ``worst-case'' disturbance
 $w_\gamma=\gamma^{-2}B_w'Xx$, where $x$ is the state of $G$, see
 \cite[Sec.\:16.8]{ZDG:95}. In other words, controller \eqref{eq:KsdHinf} generates
 the disturbance under the worst-case scenario for its analog prototype. This is
 different from the strategy proposed in \cite{LL:10}, where the sampled-data
 controller uses a piecewise-constant disturbance that ``explains'' the last deviation
 of the measured state from the calculated one. \Endoftheorem
\end{remark}

Some more observations are in order. The solvability condition of
Theorem~\ref{th:Hinf} holds for every $\gamma>\gamma_\text{opt}$ provided $\sup_ih_i$
is sufficiently small. As $\gamma\to\infty$, controller \eqref{eq:KsdHinf} recovers
the \Htwo-optimal controller of Theorem~\ref{th:H2}. If transformed to the form
of~Corollary~\ref{cor:Ksd}, controller \eqref{eq:KsdHinf} coincides with the \Hinf\ 
controller in \cite[Thm.\:5.2]{MRP:99I}, modulo replacing the sampling instances $ih$
with arbitrary $t_i$. The worst-case performance is determined by the longest sampling
interval, which is non-obvious for time-varying sampled-data systems in general.

Apropos of the worst-case sampling, the following result, whose proof is
straightforward, may be thought of as the \Hinf\ counterpart of
Proposition~\ref{pr:opttkH2}:
\begin{proposition} \label{pr:opttkHinf}
 Let $h_\gamma$ be the least upper bound for $h_i$ that satisfy the solvability
 condition of Theorem~\ref{th:Hinf} for a given $\gamma$. Then the periodic sampling
 with the sampling period $h_\gamma$ has the slowest average sampling rate among all
 sampling patterns for which the \Hinf\ performance level of $\gamma$ is attainable.
\end{proposition}


\section{Example: Design via \Hinf\ Loop Shaping} \label{sec:ex}

This section considers a numerical example, whose purpose is twofold: to illustrate
the proposed approach and to show its application to the \Hinf\ loop shaping method of
McFarlane and Glover \cite{McG:90}, which requires some light adjustments.

\subsection{Intermittent redesign for \Hinf\ loop shaping} \label{sec:HinfLSh}

The \Hinf\ loop shaping is a design procedure that uses the classical loop shaping
guidelines for choosing weights and casts the phase shaping around the crossover, the
``far from the critical point'' requirement in the classical control, as a robust
stability problem. Each iteration of this method consists of two steps.  First,
weighting functions $W_\text o$ and $W_\text i$ are chosen to shape the magnitude
(singular values) of $P_\text{msh}=W_\text oPW_\text i$. This step is technically
simple and aims at shaping loop gains in the low- and high-frequency ranges. Second, a
special robust stability problem is solved for $P_\text{msh}$ to render the
closed-loop system stable and as far from the stability margin as possible.  The
choice of the robustness setup in this step is meaningful. It is the robustness to
unstructured \Hinf\ uncertainties in the normalized coprime factors of $P_\text{msh}$.
Although normally not related to the plant physics, this problem has two important
advantages: its solution is non-iterative and it equally penalizes all four
closed-loop frequency responses (see \cite[\S4.5.1]{McG:90}). The latter means that
cancellations of stable lightly damped poles\,/\,zeros are not encouraged, in contrast
to some other optimization-based settings, like the weighted\,/\,mixed sensitivity. If
a satisfactory loop $P_\text{msh}K_0$ is reached with some choice of $W_\text o$ and
$W_\text i$ by an \Hinf\ (sub)\,optimal controller $K_0$, the resulting controller for
the original plant is $K=W_\text iK_0W_\text o$.

The robust stability problem solved in the second step is an \Hinf\ optimization
problem, whose attainable performance level may serve as a success indicator
\cite[Sec.\:6.4]{McG:90}. This renders the redesign problem of \S\ref{sec:Hinf} well
suited for this method. We actually only need to redesign $K_0$, the addition of the
weights, which are in the series connection with $K_0$, does not change the
sampled-data nature of the controller. Indeed, the series of causal and strictly
causal systems in the lifted domain is always strictly causal, see \cite[\S5.3]{MR:97}
for details.

Assume that $P_\text{msh}(s)=C(sI-A)^{-1}B$. The optimal attainable analog performance
for the \Hinf\ problem solved during the loop shaping iterations is
$\gamma_\text{opt}=\sqrt{1+\rho(YX)}$, where $X\ge0$ and $Y\ge0$ are the stabilizing
solutions to the Riccati equations (in fact, \Htwo\ Riccati equations)
\begin{gather*}
 A'X+XA+C'C-XBB'X=0, \\
 AY+YA'+BB'-YC'CY=0.
\end{gather*}
The parametrization of all $\gamma$-suboptimal solutions can then be parametrized
\cite[Thm.\:4.14]{McG:90} as $\LFTl{J_\gamma}Q$, where
\begin{equation} \label{eq:JgamHinfLSh}
 J_\gamma(s)=\mmatrix[c|cc]{A-BB'X-Z_\gamma YC'C&Z_\gamma YC'&Z_\gamma B\\\hline-B'X&0&I\\-C&I&0}
\end{equation}
and $Q$ is any linear system whose $L^2(\mathbb R^+)$-induced norm $\norm
Q<\sqrt{\gamma^2-1}$. Here $Z_\gamma\coloneq((1-\gamma^{-2})I-\gamma^{-2}YX)^{-1}>I$
is well defined for every $\gamma>\gamma_\text{opt}$. The following corollary of
Theorem~\ref{th:Hinf} can then be formulated:
\begin{corollary} \label{cor:HinfLSh}
 Let $\gamma>\sqrt{1+\rho(YX)}$. Then there is a $\gamma$-suboptimal sampled-data
 controller for a given sequence of sampling instances $\{t_i\}$ iff there exists a
 solution to the differential Riccati equation
 \begin{multline*}
  \dot P(t)=(A-YC'C)P(t)+P(t)(A'-C'CY)\\+BB'+\tfrac1{1-\gamma^{-2}}P(t)C'CP(t),\quad P(0)=Y
 \end{multline*}
 such that $\rho(P(t)X)<\gamma^2-1$, $\forall t\in[0,h_i]$ and every $i\in\mathbb
 Z^+$. If this condition holds, a sampled-data controller guaranteeing the same
 robustness level as that under $K_0$ is
 \begin{subequations} \label{eq:KsdHinfLSh}
 \begin{align}
  \dot x_\text s(t)&=Ax_\text s(t)+Bu(t)+YC'(y(t)-Cx_\text s(t)), \label{eq:KsdHinfLSh:s} \\
  \dot x_\text a(t)&=Ax_\text a(t)+Bu(t),\quad x_\text a(t_i)=Z_\gamma x_\text s(t_i) \label{eq:KsdHinfLSh:a} \\
  u(t)&=-B'Xx_\text a(t). \label{eq:KsdHinfLSh:u}
 \end{align}
 \end{subequations}
\end{corollary}
\begin{IEEEproof}
 Follows by the same steps as the proof of Theorem~\ref{th:Hinf}.
\end{IEEEproof}

Curiously, $Z_\gamma$ in \eqref{eq:KsdHinfLSh:a} is the only parameter of the
controller that depends on $\gamma$. It may be of interest to investigate the
possibility to adjust $Z_\gamma$ on-line.

\subsection{Dampening a pendulum} \label{sec:pend}

Consider the problem of controlling a pendulum, which is mounted on a cart driven by a
DC motor. The system has one input (the motor voltage) and two regulated outputs (the
cart position and the pendulum angle). Assume that the controller comprises two loops.
An internal servo loop, which is given and implemented as a 1DOF unity-feedback
system, controls the cart position. Our goal is to design the external loop, which
aims at dampening pendulum oscillations during command response of the cart. The
external loop measures the pendulum angle and modifies the reference signal to the
inner loop. This way, the reference signal for the cart is treated as the load
disturbance against which the external loop acts.

\begin{figure}[!tb]
 \centering
  \subfigure[Pendulum angle (dashed curve represents the open-loop response)]{\label{fig:ya}\includegraphics[scale=0.667]{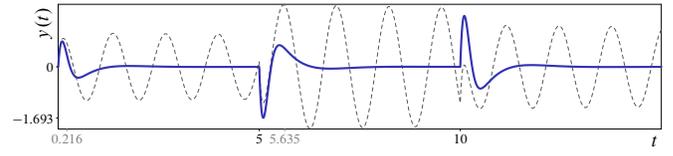}} \\
  \subfigure[Control input for the pendulum loop (corrections to the reference)]{\label{fig:ua}\includegraphics[scale=0.667]{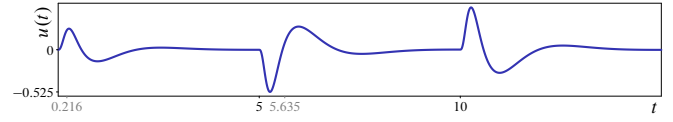}}
 \caption{Responses to a square wave, analog $K_0$ designed for $\gamma=3.703$} \label{fig:yua}
\end{figure}

Let the transfer function from the servo reference signal to the pendulum angle be
\[
 P(s)=-\frac{42s^2}{(s+18)(s^2+0.02s+23)}.
\]
It has a pair of lightly damped poles at $s=-0.01\pm\J4.796$, so the control goal is
to dampen them by feedback. To this end, we design an analog controller via the \Hinf\ 
loop shaping procedure. The choice
\[
 W_\text i(s)=\frac5{s+2}
 \quad\text{and}\quad
 W_\text o(s)=1
\]
yields a satisfactory loop with low $\gamma_\text{opt}=1.7213$. Consider then the
design with $\gamma=3.703\approx2.151\gamma_\text{opt}$ (the rationale behind this
choice will be clarified later on), which produces the central analog controller
\[
 K_0(s)=W_\text i(s)\frac{12.534(s+18.85)(s+1.839)(s+0.2895)}{(s^2+1.91s+1.514)(s^2+37.26s+547.4)}.
\]
The response of the resulted closed-loop system to a square wave load disturbance with
a magnitude of $\pm0.5$ and a period of $10$\,sec, is shown in Fig.\,\ref{fig:yua} by
solid blue lines.
Dampening properties the designed feedback are apparent from
comparing the closed-loop output response to that of the open-loop plant (dashed line
in Fig.\,\ref{fig:ya}).

\begin{figure}[!h]
 \centering\includegraphics[scale=0.67]{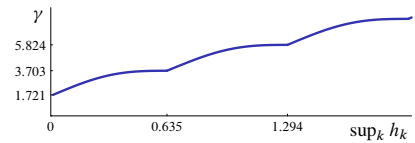}
 \caption{Attainable $\gamma$ as a function of the largest sampling interval} \label{fig:gammavstk}
\end{figure}

To redesign $K_0$, consider first how the condition of Corollary~\ref{cor:HinfLSh} on
$\{t_i\}$ depend on the robustness level $\gamma$. Calculating the least upper bound
on the admissible sampling period at each $\gamma>\gamma_\text{opt}$, we end up with
the plot in Fig.\,\ref{fig:gammavstk}.
Expectably, the required $\sup_ih_i$ for $\gamma$'s close to $\gamma_\text{opt}$ is
quite close to zero, which leaves little room for investigating properties of
intermittent sampling. It therefore makes sense to consider larger $\gamma$. The value
chosen in the design of $K_0$ is at the point where the slope of the curve in
Fig.\,\ref{fig:gammavstk} is zero (so minimal damage for the increase of $h_i$). The
maximal admissible sampling period in this case is $0.635$, which is rather slow from
the classical sampled-data control viewpoint, as the corresponding Nyquist frequency
of almost $5$\,rad/sec is comparable with the largest loop crossover of
$7.75$\,rad/sec, see also the transients in Fig.\,\ref{fig:yua}.

Having the bound for admissible sampling rates and complete freedom in choosing the
sampling pattern within this bound, let us dream up the following strategy for the
choice of $t_i$. Consider the signal $\eta=Q_\text{stat}(y-Cx_\text s)$, where
$Q_\text{stat}$ is given by \eqref{eq:QstatHinf}, adopted to $J_\gamma$ in
\eqref{eq:JgamHinfLSh}. This signal is reset at every sampling instance $t_i$. As the
norm of this $Q_\text{stat}$ determines the \Hinf\ performance, we may use the
$L^2$-norm of $\eta$ as a basis for event generation. To this end, let $\theta_i$ be
the solution of
\[
 \int_0^\theta\eta'(t_i+t)\eta(t_i+t)\D t=0.025^2
\]
and consider the following sampling generation mechanism:
\[
 h_i=\min\bigl\{\,\theta_i,0.635\,\bigr\},
\]
which is easy to implement. In other words, the controller samples either as the $L^2$
norm of $\eta$ reaches $0.025$ or after $0.635$\,sec if the norm does not reach this
level by then.

\begin{figure}[!tb]
 \centering
  \subfigure[Pendulum angle]{\label{fig:ysd}\includegraphics[scale=0.667]{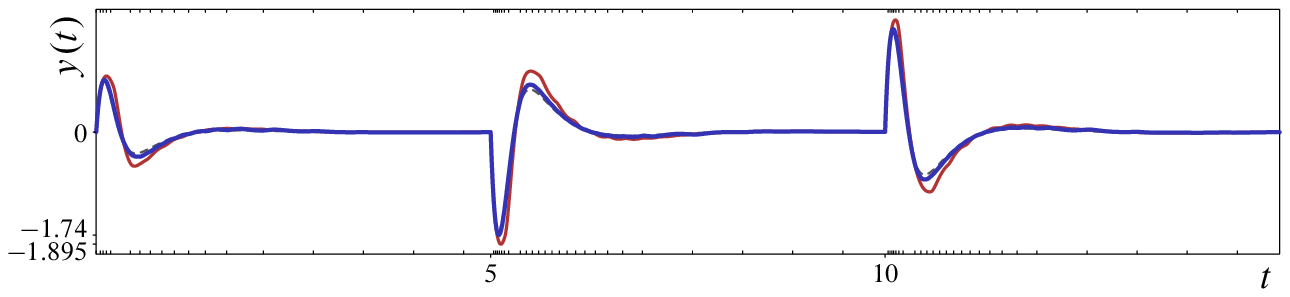}} \\
  \subfigure[Control input for the pendulum loop (corrections to the reference)]{\label{fig:usd}\includegraphics[scale=0.667]{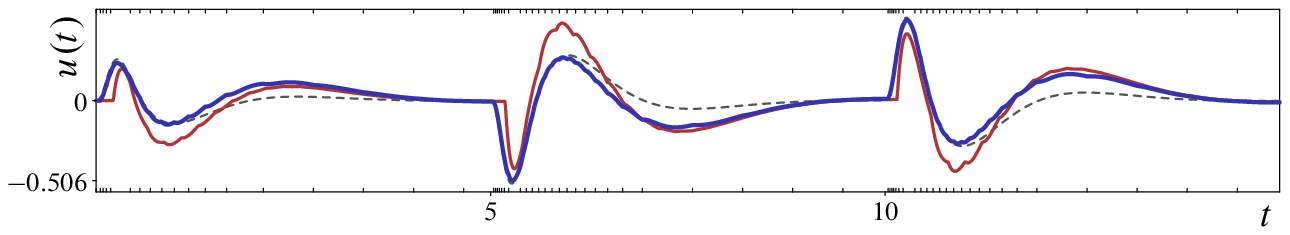}}
 \caption{Responses to a square wave, intermittent redesign of $K_0$ (blue lines: event-based sampling, marked as the $x$-axis ticks; red lines: uniform sampling with the same density; gray dashed lines: analog controller)} \label{fig:yusd}
\end{figure}

Simulation results with this controller are presented in Fig.\,\ref{fig:yusd} by blue
lines.
The resulted sampling instances are marked as the $x$-axis ticks. Intuitively, the
sampling rate increases during the transients and decreases as the steady state is
reached. One can see that the output response is quite close to the response under the
analog $K_0$ (dashed gray line in Fig.\,\ref{fig:ysd}). This is noteworthy, taking
into account that the average sampling period here, $h_\text{av}=0.216$, is still
rather slow (the corresponding Nyquist frequency, $14.5$\,rad/sec, exceeds the largest
crossover of the analog loop only by a factor of $2$). For the sake of comparison, the
red lines in Fig.\,\ref{fig:yusd} present responses under periodic sampling with
$h_i=0.216$, $\forall i$. Note that the control signals $u(t)$ are continuous
functions under both sampling strategies. This is because the discontinuous signal
generated by \eqref{eq:KsdHinfLSh:u} is then filtered by the low-pass $W_\text i$. A
larger pole excess in $W_\text i(s)$ would result in a differentiable $u(t)$.


\section{Concluding Remarks} \label{sec:concl}

The paper has studied the problem of digital redesign of analog controllers under
intermittent, possibly unknown a priori, sampling. The main idea, borrowed from
\cite{MR:97}, is to use the characterization of causal sampled-data controllers as the
set of all strictly causal systems in the lifted domain to extract sampled-data
controllers from Youla-like parametrizations of satisfactory analog controllers. The
resulting controllers are always stabilizing and, if optimal control parametrizations
are considered, performance guaranteeing. As a byproduct of the proposed approach, the
\Htwo\ and \Hinf\ problems under intermittent sampling have been solved. In both cases
the (sub)\,optimal control laws are explicit and readily computable. It has also been
proved that the uniform sampling is optimal among all sampling patterns with a given
sampling density.

Some extensions of the results put forward in this paper should be immediate. For
example, adding a single loop delay can be addressed via the loop shifting approach,
similarly to the treatment of the constant sampling rate in \cite{MST:14}. This way
both stabilization and \Htwo\ optimization problems can be solved, thus justifying the
predictor-based structure proposed in \cite{GLLG:11} without a proof. This approach
will not work in the \Hinf\ case though. Another alternation that seems to be
immediate is to apply the ideas of this paper to the formulation proposed in
\cite[Ch.\:4]{GAM:14}, where the analog loop is closed not only instantaneously, but
rather during some short time intervals. A more laborious extension would be to come
up with a theoretically justified event generation mechanism.


\appendix

\subsection{Proof of Theorem~\ref{th:H2}} \label{sec:H2proof}

We start with the following technical result:
\begin{lemma} \label{l:HtwoCLN}
 Let $J_0$ be given by \eqref{eq:J0ss:ob} with $F$ and $L$ as in the statement of
 Theorem~\eqref{th:H2}. Consider the family of controllers $\LFTl{J_0}Q$ for a causal
 linear $Q$ such that $\norm Q_2<\infty$. Then
 \[
  \norm{T_{zw}}_2^2=\gamma_0^2+\norm Q_2^2,
 \]
 where $\gamma_0$ is the optimal \Htwo\ performance attainable by contin\-uous-time
 controllers.
\end{lemma}
\begin{IEEEproof}
 The closed-loop map for the considered family of controllers is
 \cite[Thm.\:12.16]{ZDG:95} $T_{zw}=T_1+T_2QT_3$, where
 \[
  \mmatrix{T_1(s)&T_2(s)\\T_3(s)&0}=\mmatrix[cc|cc]{A_F&-B_uF&B_w&B_u\\0&A_L&B_L&0\\\hline C_F&-D_{zu}F&0&D_{zu}\\0&C_y&D_{yw}&0}
 \]
 with Hurwitz $A_F\coloneq A+B_uF$ and $A_L\coloneq A+LC_y$, $B_L\coloneq
 B_w+LD_{yw}$, and $C_F\coloneq C_z+D_{zu}F$. Moreover, $T_1\in\Htwo$, $T_2$ is inner
 \cite[Thm.\:13.32]{ZDG:95} and $T_3$ is co-inner \cite[Thm.\:13.35]{ZDG:95}.

 Now, \eqref{eq:H2norm} defines a (degenerate) Hilbert space with the inner product
 \[
  \inpr{G_1}{G_2}_2=\lim_{T\to\infty}\frac1T\int_0^T\int_\tau^\infty\trace(g_2'(t,\tau)g_1(t,\tau))\,\D t\D\tau,
 \]
 so that $\norm{G}_2^2=\inpr GG_2$. If $G$ is a causal LTI system, its adjoint with
 respect to the inner product above, $G^*$, is the anti-causal LTI system, whose
 transfer function equals $[G(-s)]'$, exactly as in the case of the conventional
 \Htwo\ space. We then have:
 \begin{align*}
  \norm{T_{zw}}_2^2
   &=\inpr{T_1+T_2QT_3}{T_1+T_2QT_3}_2 \\
   &=\norm{T_1}_2^2+\norm{T_2QT_3}_2^2+2\real\inpr{T_2QT_3}{T_1}_2 \\
   &=\norm{T_1}_2^2+\norm{Q}_2^2+2\real\inpr{Q}{V}_2,
 \end{align*}
 where $V\coloneq T_2^*T_1T_3^*$ and the facts that $T_2^*T_2=I$ and $T_3T_3^*=I$ were
 used. It can be verified, via straightforward state-space manipulations, that $V$ is
 anti-causal, with
 \[
  V(s)=\mmatrix[cc|c]{-A_F'&-A_F'XY-XAY&XL\\0&-A_L'&C_y'\\\hline B_u'&-D_{zu}'C_zY&0},
 \]
 where $X\ge0$ and $Y\ge0$ are the stabilizing solutions of the state-feedback and
 filtering Riccati equations, respectively. This implies that the responses of $V$ and
 $Q$ to the same impulse have disjoint supports. Therefore, $\inpr{Q}{V}_2=0$, which
 completes the proof (with $\gamma_0=\norm{T_1}_2$).
\end{IEEEproof}

By Lemma~\ref{l:stabsdld}, the controller of the form $\LFTl{J_0}Q$ is a sampled-data
one iff $Q=Q_\text{stat}+Q_\text{sd}$ for a given $Q_\text{stat}$ and any stable
sampled-data $Q_\text{sd}$. Remember that the lifting of $Q_\text{stat}$ is static and
the lifting of $Q_\text{sd}$ is strictly proper. Therefore, the impulse responses of
$Q_\text{stat}$ and $Q_\text{sd}$ are non-overlapping for any admissible
$Q_\text{sd}$, which, in turn, implies that
\[
 \norm Q_2^2=\norm{Q_\text{stat}+Q_\text{sd}}_2^2=\norm{Q_\text{stat}}_2^2+\norm{Q_\text{sd}}_2^2.
\]
Thus, the optimal performance is attained with any $Q$ such that $Q-Q_\text{stat}$ is
in the kernel of semi-norm \eqref{eq:H2norm}.

Compute now $\norm{Q_\text{stat}}_2^2$. By \eqref{eq:Qstatss}, $Q_\text{stat}$ can be
described by
\begin{equation} \label{eq:QstatH2}
 \begin{aligned}
  \dot x_Q(t)&=Ax_Q(t)+L\epsilon(t),\quad x_Q(t_i)=0 \\
  \eta_Q(t)&=Fx_Q(t)
 \end{aligned}
\end{equation}
Its impulse response is
$q_\text{stat}(t,\tau)=F\E^{A(t-\tau)}L\mathbb1_{[\tau,t_j)}(t)$, where $t_j$ is the
smallest element of $\{t_i\}$ such that $t_j\ge\tau$ and $\mathbb1_{[a,b)}(t)$ is the
characteristic function of the interval $[a,b)$. Then
\begin{align*}
 \norm{Q_\text{stat}}_2^2
  &=\lim_{i\to\infty}\frac1{t_i}\int_0^{t_i}\int_\tau^\infty\normf{q_\text{stat}(t,\tau)}^2\,\D t\D\tau \\
  &=\lim_{i\to\infty}\frac1{t_i}\sum_{j=0}^{i-1}\int_{t_j}^{t_{j+1}}\int_\tau^\infty\normf{q_\text{stat}(t,\tau)}^2\,\D t\D\tau \\
  &=\lim_{i\to\infty}\frac1{t_i}\sum_{j=0}^{i-1}\int_{t_j}^{t_{j+1}}\int_\tau^{t_{j+1}}\normf{F\E^{A(t-\tau)}L}^2\,\D t\D\tau,
\end{align*}
from which the expression for the achievable performance follows by straightforward
integration variable change.

Finally, the optimal control law is in form \eqref{eq:Ksdob} because $K_0$ is observer
based. \hfill\IEEEQEDclosed

\subsection{Proof of Theorem~\ref{th:Hinf}} \label{sec:Hinfproof}

In addition to the notation introduced prior to the formulation of the Theorem, define
\[
 \tilde B_u\coloneq B_u+\gamma^{-2}YC_z'D_{zu},\quad \tilde C_y\coloneq C_y+\gamma^{-2}D_{yw}B_w'X,
\]
and $Z_\gamma\coloneq(I-\gamma^{-2}YX)^{-1}$. It is known \cite[Thm.\:16.5]{ZDG:95}
that if $\gamma>\gamma_\text{opt}$, all $\gamma$-suboptimal LTI controllers can be
characterized as $\LFTl{J_\gamma}Q$ for
\begin{equation} \label{eq:J0Hinf}
 J_\gamma(s)=\mmatrix[c|cc]{A_\gamma&-Z_\gamma L&Z_\gamma\tilde B_u\\\hline F&0&I\\-\tilde C_y&I&0}
\end{equation}
and an arbitrary LTI $Q\in\Hinf$ such that $\norm Q_\infty<\gamma$, where
$A_\gamma\coloneq A+\gamma^{-2}B_wB_w'X+B_uF+Z_\gamma L\tilde C_y$. Because the
central controller is the one corresponding to $Q=0$, $K_0=\LFTl{J_\gamma}0$. The
parametrization above extends to time-varying controllers as well. Namely, the set of
all $\gamma$-suboptimal linear causal controllers is $\LFTl{J_\gamma}Q$, where $Q$ is
an arbitrary bounded causal operator on $L^2(\mathbb R^+)$ such that its induced norm
$\norm Q<\gamma$, see the arguments in \cite{Tad:90}.

By Lemma~\ref{l:stabsdld}, a controller of the form $\LFTl{J_\gamma}Q$ is in the
sampled-data form iff $Q=Q_\text{stat}+Q_\text{sd}$ for a $Q_\text{stat}$, verifying
\begin{equation} \label{eq:QstatHinf}
 \begin{aligned}
  \dot x_Q(t)&=A_\gamma^\times x_Q(t)+Z_\gamma L\epsilon(t),\quad x_Q(t_i)=0 \\
  \eta(t)&=Fx_Q(t)
 \end{aligned}
\end{equation}
where $A_\gamma^\times\coloneq A_\gamma-Z_\gamma(\tilde B_uF+L\tilde
C_y)=A+\gamma^{-2}(B_wB_w'X+Z_\gamma YF'F)$, and any stable causal sampled-data
$Q_\text{sd}$. The existence of an admissible $Q$ is then equivalent to the existence
of a causal sampled-data system $Q_\text{sd}$ such that
$\norm{Q_\text{stat}+Q_\text{sd}}<\gamma$. To address the latter, the following result
is required:
\begin{lemma} \label{l:HinfQ}
 $\norm{Q_\text{stat}+Q_\text{sd}}\ge\norm{Q_\text{stat}}$ for all causal sampled-data
 systems $Q_\text{sd}$.
\end{lemma}
\begin{IEEEproof}
 In the lifted domain, $\vl Q_\text{stat}$ is static and $\vl Q_\text{sd}$ is strictly
 causal. Hence, the responses of $\vl Q_\text{stat}$ and $\vl Q_\text{sd}$ to any
 input $\vl\epsilon$ such that $\vl\epsilon[i]=0$ for all $i\ne j$ for some given
 $j\in\mathbb Z^+$ are non-overlapping (zeros $\forall i\ne j$ and $\forall i\le j$,
 respectively). As a result, in the time domain we have that for any $\epsilon(t)$
 with support in $[t_i,t_{i+1})$,
 \[
  \norm{(Q_\text{stat}+Q_\text{sd})\epsilon}_2^2=\norm{Q_\text{stat}\epsilon}_2^2+\norm{Q_\text{sd}\epsilon}_2^2\ge\norm{Q_\text{stat}\epsilon}_2^2.
 \]
 where $\norm\cdot_2$ stands for the $L^2(\mathbb R^+)$ signal norm. The result then
 follows by observing that the worst-case input for $Q_\text{stat}$ has support in
 $[t_i,t_{i+1})$ for some $i$, which, in turn, is a consequence of the fact that
 $Q_\text{stat}$ resets at each $t_i$ (by Lemma~\ref{l:spls}).
\end{IEEEproof}

It follows from Lemma~\ref{l:HinfQ} that an admissible $Q$ exists iff
$\norm{Q_\text{stat}}<\gamma$ (as we can always pick $Q_\text{sd}=0$). The norm bound
can then be verified by the following result:
\begin{lemma} \label{l:normQstat}
 Let $\gamma>\gamma_\text{opt}$ and $Q_\text{stat}$ be given by \eqref{eq:QstatHinf}.
 Then $\norm{Q_\text{stat}}<\gamma$ iff the conditions of the Theorem hold.
\end{lemma}
\begin{IEEEproof}
 It is readily seen that $\norm{Q_\text{stat}}<\gamma$ iff the $L^2[0,h_i)$-induced
 norm of $\LFTu{J_\gamma^{-1}(s)}0=F(sI-A_\gamma^\times)^{-1}Z_\gamma L$ is less than
 $\gamma$ for all $i\in\mathbb Z^+$. But the $L^2[0,h)$-induced norm of an LTI system
 is a monotonically increasing function of $h$. Hence, we only need to check the norm
 for the maximal $h_i$.
 
 It is known \cite[Lem.\:2.2]{GCT:96} that the $L^2[0,h)$-induced norm of
 $\LFTu{J_\gamma^{-1}}0$ is less than $\gamma$ iff the differential Riccati equation
 \[
  \dot R(t)=A_\gamma^\times R(t)+R(t)(A_\gamma^\times)'+Z_\gamma LL'Z_\gamma'+\gamma^{-2}R(t)F'FR(t)
 \]
 with $R(0)=0$ has a bounded solution in the whole interval $[0,h]$. This Riccati
 equation, in turn, is associated with the Hamiltonian matrix \cite[Lem.\:2.3]{GCT:96}
 \[
  H_R\coloneq\mmatrix{-(A_\gamma^\times)'&-\gamma^{-2}F'F\\Z_\gamma LL'Z_\gamma'&A_\gamma^\times}.
 \]
 It can be shown \cite[Eqn.\:(14)]{Mir:03tac} that
 \[
  H_R=\mmatrix{Z_\gamma'&\gamma^{-2}X\\YZ_\gamma'&I}^{-1}H_P\mmatrix{Z_\gamma'&\gamma^{-2}X\\YZ_\gamma'&I},
 \]
 where
 \[
  H_P\coloneq\mmatrix{-A'&-\gamma^{-2}C_z'C_z\\B_wB_w'&A}
 \]
 is the Hamiltonian matrix associated with $P(t)$. As a result,
 \[
  R(t)=(I-\gamma^{-2}P(t)X)^{-1}(P(t)-Y)Z_\gamma',
 \]
 so that it is bounded iff $\det(I-\gamma^{-2}P(t)X)\ne0$. It is readily seen that
 $P_\text d(t)\coloneq\dot P(t)$ satisfies the Lyapunov equation
 \[
  \dot P_\text d(t)=A_P(t)P_\text d(t)+P_\text d(t)A_P'(t),\quad P_\text d(0)=LL'\ge0
 \]
 where $A_P\coloneq A+\gamma^{-2}PC_z'C_z$. Hence, $\dot P(t)\ge0$ for all $t$ and
 $P(t)$ is non-decreasing. We also know that $\rho(P(0)X)<\gamma^2$ whenever
 $\gamma>\gamma_\text{opt}$. Thus, the boundedness of $R(t)$ in $[0,h]$ is equivalent
 to $\rho(P(t)X)<\gamma^2$ at each $t$ in this interval.
\end{IEEEproof}

To complete the proof of the Theorem, we only need to show that controller
\eqref{eq:KsdHinf} is a particular case of \eqref{eq:Ksd} if $J_0=J_\gamma$. This can
be verified by direct substitution using the fact that
\[
 A_\gamma-Z_\gamma\tilde B_uF=Z_\gamma A_LZ_\gamma^{-1},
\]
which can be verified via some lengthly algebra. \hfill\IEEEQEDclosed


\section*{Acknowledgments}

I am indebted to Igor Gindin and Miriam Zacksenhouse for drawing my attention to the
problem via \cite{GLLG:11}.


\end{document}